\documentclass[final]{siamltex}
\usepackage{lmodern,eurosym,ae,tabularx,amssymb,url,epsf,epsfig,amsmath,amsfonts,graphicx,dsfont}
\font\dsrom=dsrom10 scaled 1200           
\newtheorem{thm}{Theorem}[section]
\newtheorem{prop}[thm]{Proposition}
\newtheorem{lemme}[thm]{Lemma}

\newtheorem{rmq}[thm]{Remark}
\def\strike#1{{\setbox0=\hbox{\kern3pt #1\kern3pt}\copy0\kern-\wd0
\setbox1=\hbox{\raise2.5pt\hbox{\vrule height.4pt depth0pt width\wd0}}\copy1}}
\newenvironment{pf}{\ \\ {\it Proof}}{\hfill\mbox{$\diamond$}\medskip}
\renewcommand{\P}{\mathbb{P}}
\newcommand{\indicatrice}{\textrm{\dsrom{1}}} 
\newcommand{\F}{{\mathcal F}}
\newcommand{\R}{{\mathbb R}}
\newcommand{\N}{{\mathbb N}}
\newcommand{\Z}{{\mathbb Z}}
\newcommand{\E}{{\mathbb{E}}}
\renewcommand{\P}{{\mathbb{P}}}

\title{Continuity correction for barrier options in jump-diffusion models}
\author{El Hadj Aly Dia\thanks{Universit\'e Paris-Est, Laboratoire d'Analyse et de Math\'ematiques Appliqu\'ees, UMR CNRS 8050, 5 bd. Descartes, Champs-sur-Marne, $77454$ Marne-la-Vall\'ee, France ({\tt dia.eha@gmail.com}).} \and Damien Lamberton\thanks{Universit\'e Paris-Est, Laboratoire d'Analyse et de Math\'ematiques Appliqu\'ees, UMR CNRS 8050, 5 bd. Descartes, Champs-sur-Marne, $77454$ Marne-la-Vall\'ee, France ({\tt damien.lamberton@univ-mlv.fr}).}}
\date{}
\begin{document}

\maketitle
  
\begin{abstract}
The aim of this paper is to study the continuity correction for barrier options in jump-diffusion models. 
For this purpose, we express the pay-off of a barrier option in terms of the maximum
of the underlying process. We then condition with respect to the jump times and to the values
of the underlying at the jump times and rely on the connection between 
the maximum of the Brownian motion and Bessel processes.
\end{abstract}
\begin{keywords} 
Barrier option, Bessel process, Continuity correction, Exponential L\'evy model, jump-diffusion.
\end{keywords}
\begin{AMS}
60G51, 60J75, 91G20
\end{AMS}
\begin{JEL}
C02, G13
\end{JEL}
\pagestyle{myheadings}
\thispagestyle{plain}
\markboth{E. H. A. DIA AND D. LAMBERTON}{CONTINUITY CORRECTION FOR BARRIER OPTIONS}
\section{Introduction}
\label{sec:intro}

 In the Black-Scholes setting,  Broadie, Glasserman and Kou (1997) and Kou (2003)
  derived continuity correction 
 formulas for  barrier options. The purpose of this paper is to establish similar results for 
 jump-diffusion models. 
 The approach of Broadie, Glasserman and Kou was based on the expression of the pay-off of a
barrier option in terms of the hitting time of the barrier by the underlying stock price. 
They managed to relate
the discrete barrier option price to the continuous one by using classical results on the {\em overshoot}
asymptotics of the Gaussian random walk.

Our approach is completely different and provides a new proof of the Broadie-Glassserman-Kou
results, even in the Black-Scholes case. We start from the expression of the pay-off of barrier options
 in terms of the 
maximum process, which essentially involves the cumulative distribution function of the maximum. 
We then rely on the connection between the maximum of Brownian motion and the Bessel process,
 following the ideas of Asmussen, Glynn, Pitman (1995), in their study of the weak convergence
 of the normalized difference between the continuous and discrete maximum of Brownian motion.
 The extension to jump-diffusions is obtained by conditioning with respect to the jump times
 and to the values of the process at the jump times.

Note that the Asmussen-Glynn-Pitman Theorem was the basic tool for the derivation 
by Broadie, Glasserman and Kou (1999) of continuity corrections
for {\em lookback} options, and we showed in \cite{DiaLamberton} that this approach could
be extended to jump-diffusion processes. 
The dependence of the payoff with respect to the maximum
is much less smooth in the case of barriers, and we will need to go deeper into the connection
between the maximum and the Bessel process to prove our results. In some sense, our results
prove that continuity correction formulas can be obtained in a unified way for barrier and for lookback
options.

The paper is organized as follows. In the next section, we 
present our main results: a continuity correction formula for a general pay-off (see 
Theorem~\ref{probcorrection}), and its application to barrier options (see
Proposition~\ref{barriercorrection}). We also demonstrate the use of these results
by showing some numerical results for a double-exponential jump-diffusion model.  
The other sections of the paper 	are devoted to the proof of Theorem~\ref{probcorrection}.
 In Section~\ref{sec:Poisson},
we derive some preliminary estimates on the jump times of a Poisson process.
In Section~\ref{sec:jump}, by conditioning with respect to the jump times, we reduce the problem
to the analysis of discrete vs continuous supremum between the jump times. 

In Section~\ref{section-rep},
we further condition with respect to the values of the underlying process at the jump times.
We then have to deal with independent Brownian motions, and we establish a representation
of a conditional expectation of a function of the maximum, the discrete maximum and the terminal value
in terms of Bessel processes (see Proposition~\ref{prop-rep-Bessel}).

Section~\ref{sec:Bessel} is devoted to the derivation of some elementary estimates concerning
the transition kernel of the Bessel process which are needed in the last two sections.
In Section~\ref{sec:domination}, we derive some bounds 
for conditional expectations, in order to be able to derive convergence results 
for the unconditional expectations from the corresponding results for conditional expectations.
In Section~\ref{sec:convergence}, we establish the continuity correction for conditional expectations.
\section{Continuity correction formulas}
\label{sec:prelim}
In a jump-diffusion model, the price of the underlying stock at time $t$ is given by
\begin{eqnarray*}
S_t=S_0e^{X_t},\quad 0\leq t\leq T,
\end{eqnarray*}
where, under the pricing measure, the process $X=(X_t)_{0\leq t\leq T}$ is given by
\begin{equation}\label{levyactivitefinie}
X_t=\gamma t+\sigma B_t + \sum_{i=1}^{N_t}Y_i,
\end{equation}
where $\gamma$ and $\sigma$ are real constants, with $\sigma>0$,
$(B_t)_{0\leq t\leq T}$ is a standard Brownian motion,
$N$ is a Poisson process with intensity $\lambda>0$, 
and $\left(Y_i\right)_{i\geq1}$ are i.i.d. random variables.
Note that, under the pricing measure, the process $(e^{-(r-\delta)t}S_t)_{0\leq t\leq T}$,
where $r$ is the interest rate and $\delta$ the dividend rate
is a martingale. This implies the following relation between $\gamma$ and the other parameters
\[
\gamma=r-\delta-\frac{\sigma^2}{2}+\lambda \E\left(e^{Y_1}-1\right).
\]
In the terminology of exponential L\'evy models, note that $X$ is a L\'evy process
with a non-zero Brownian part and a finite L\'evy measure, given by $\nu(dx)=\lambda \mu(dx)$,
where $\mu$ is the distribution of the random variable $Y_1$.
 For more details about L\'evy processes see \cite{sato}.

We define
\begin{eqnarray*}
 &&M_t^{X}=\sup_{0\leq s\leq t}X_{s}, \ \ M_t^{X,n}=\max_{0\leq k\leq n}X_{\frac{kt}{n}}
 \\&&m_t^X=\inf_{0\leq s\leq T}X_{s}, \ m_t^{X,n}=\min_{0\leq k\leq n}X_{\frac{kt}{n}}.
\end{eqnarray*}
When there is no ambiguity we can remove the super index $X$. 


The options we will consider in the sequel will have as underlying the asset with price $S$. 
We will denote by $K$ and $H$ the strike and the barrier of the option. 
The maturity of the options is assumed to be $T$.
Figures \ref{payoffscall_tab} and \ref{payoffsput_tab} give the payoffs of barrier options. The corresponding prices are the expected values of the discounted payoffs.

\begin{figure}[ht]
\begin {center}
\begin{tabular}{|c|c|c|}
 \hline 
 Barrier &Continuous&Discrete\\ \hline
 \emph{Up Out}&$\left(S_T-K\right)^{+}\indicatrice_{\left\{S_0e^{M_T}< H\right\}}$&$\left(S_T-K\right)^{+}\indicatrice_{\left\{S_0e^{M_T^n}< H\right\}}$\\ \hline
 \emph{Up In}&$\left(S_T-K\right)^{+}\indicatrice_{\left\{S_0e^{M_T}\geq H\right\}}$&$\left(S_T-K\right)^{+}\indicatrice_{\left\{S_0e^{M_T^n}\geq H\right\}}$\\ \hline
\emph{Down Out}&$\left(S_T-K\right)^{+}\indicatrice_{\left\{S_0e^{m_T}> H\right\}}$&$\left(S_T-K\right)^{+}\indicatrice_{\left\{S_0e^{m_T^n}> H\right\}}$\\ \hline
\emph{Down In}&$\left(S_T-K\right)^{+}\indicatrice_{\left\{S_0e^{m_T}\leq H\right\}}$& $\left(S_T-K\right)^{+}\indicatrice_{\left\{S_0e^{m_T^n}\leq H\right\}}$\\ \hline
\end{tabular}
\caption{Payoffs of barrier call options}
\label{payoffscall_tab}
\end {center}
\end{figure}

\begin{figure}[ht]
\begin {center}
\begin{tabular}{|c|c|c|}
 \hline 
 Barrier &Continuous&Discrete\\ \hline
 \emph{Up Out}&$\left(K-S_T\right)^{+}\indicatrice_{\left\{S_0e^{M_T}< H\right\}}$&$\left(K-S_T\right)^{+}\indicatrice_{\left\{S_0e^{M_T^n}<
  H\right\}}$\\ \hline
\emph{Up In}& $\left(K-S_T\right)^{+}\indicatrice_{\left\{S_0e^{M_T}\geq H\right\}}$&$\left(K-S_T\right)^{+}\indicatrice_{\left\{S_0e^{M_T^n}\geq H\right\}}$\\ \hline
\emph{Down Out}& $\left(K-S_T\right)^{+}\indicatrice_{\left\{S_0e^{m_T}> H\right\}}$&$\left(K-S_T\right)^{+}\indicatrice_{\left\{S_0e^{m_T^n}> H\right\}}$\\ \hline
\emph{Down In}& $\left(K-S_T\right)^{+}\indicatrice_{\left\{S_0e^{m_T}\leq H\right\}}$&$\left(K-S_T\right)^{+}\indicatrice_{\left\{S_0e^{m_T^n}\leq H\right\}}$\\ \hline
\end{tabular}
\caption{Payoffs of barrier put options}
\label{payoffsput_tab}
\end {center}
\end{figure}
Let $UOC(H)$ be the price of a continuous \emph{up and out} call  with barrier $H$, We have
\begin{eqnarray*}
UOC(H)=\E e^{-rT}\left(S_0e^{X_T}-K\right)^{+}\indicatrice_{\left\{\sup_{0\leq t\leq T}S_0e^{X_t}<H\right\}}.
\end{eqnarray*}
Define
\[
k=\log\left(\frac{K}{S_0}\right),\quad h=\log\left(\frac{H}{S_0}\right).
\]
We can write
\begin{eqnarray*}
UOC(H)&=&S_0\E e^{-rT}e^{X_T}\indicatrice_{\left\{M_T<h,X_T>k\right\}}-Ke^{-rT}\P\left[M_T<h,X_T>k\right]
\\&=&S_0e^{-\delta T}\E e^{-(r-\delta)T}e^{X_T}\indicatrice_{\left\{M_T<h,X_T>k\right\}}-Ke^{-rT}\P\left[M_T<h,X_T>k\right].
\end{eqnarray*}
We know that the process $\left(e^{-(r-\delta)t}e^{X_t}\right)_{0\leq t\leq T}$ is a martingale. 
Let $\bar{\P}$ be the probability defined by its density with respect to the 
pricing probability measure  $\P$ 
\begin{eqnarray*}
\frac{d \bar{\P}}{d\P}=e^{-(r-\delta)T +X_T}.
\end{eqnarray*}
Note that (as can be deduced, for instance, from  Theorem $3.9$ of \cite{kyprianou}),
the process $X$ remains a L\'evy process under probability $\bar{\P}$, and that its L\'evy measure
under $\bar{\P}$ is given by $\bar{\nu}(dx)=e^{x}\nu(dx)$.

We have
\begin{eqnarray*}
UOC(H)&=&S_0e^{-\delta T}\bar{\P}\left[M_T<h,X_T>k\right]-Ke^{-rT}\P\left[M_T<h,X_T>k\right]
\end{eqnarray*}
If we call $UOC^n$ the price of a discrete \emph{up and out} call  with barrier $H$, and $n$ fixing dates (with step $\frac{T}{n}$), then we have similarly
\begin{eqnarray*}
UOC^n(H)&=&S_0e^{-\delta T}\bar{\P}\left[M_T^n<h,X_T>k\right]-Ke^{-rT}\P\left[M_T^n<h,X_T>k\right].
\end{eqnarray*}
Finding continuity corrections between continuous and discrete barrier options
amounts in fact to finding corrections between the above probabilities. This is the aim of the following result.

\begin{thm}\label{probcorrection}
Let $X$ be an integrable  L\'evy process of the form (\ref{levyactivitefinie}), with $\sigma>0$. For
any bounded Borel measurable function $g:\R\to \R$ and for any positive number $x$,
we have
\begin{eqnarray*}
\E\left[g(X_T)\indicatrice_{\{M_T\geq x>M_T^n\}}\right]=
\E\left[g(X_T)\indicatrice_{\{M_T\geq x> M_T-\frac{\sigma\sqrt{T}\beta_1}{\sqrt{n}}\}}\right]+
  o\left(\frac{1}{\sqrt{n}}\right),
\end{eqnarray*}
where $\beta_1=\E R$ and $R$ is defined by
\begin{equation}\label{R}
R=\min_{\{j\in\Z\}}\check{R}(U+j).
\end{equation}
Here, $(\check {R}(t))_{t\in \R}$ is a two sided three dimensional Bessel process (i.e. $\check {R}(t)=R_1(t)$ for $t\geq 0$ and $\check R(t)=R_2(-t)$ for $t<0$, where $R_1$ and $R_2$ are independent copies of the usual three dimensional Bessel process, starting from $0$) and $U$ is uniformly distributed on $[0,1]$ and independent of $\check{R}$. 
\end{thm}

Note that the result does not depend on the jump part of the process, so that the continuity correction
for jump-diffusion models is the same as for the Black-Scholes model.

The result of Theorem~\ref{probcorrection} can also be written in the form
\begin{eqnarray*}
\E\left[g(X_T)\indicatrice_{\{M_T< x+\frac{\sigma\sqrt{T}\beta_1}{\sqrt{n}}\}}\right]=
\E\left[g(X_T)\indicatrice_{\{M_T^n< x\}}\right]+o\left(\frac{1}{\sqrt{n}}\right).
\end{eqnarray*}
Moreover, the proof of Theorem~\ref{probcorrection} shows that the theorem 
is still true if we replace $x$ by a sequence $x_n$ which converges to $x$ when $n\rightarrow+\infty$.
 So, under the assumptions of Theorem~\ref{probcorrection}, we have
\begin{eqnarray}
\P\left(\!M_T<x+\frac{\sigma\sqrt{T}\beta_1}{\sqrt{n}},X_T>y\!\right)&=&\P\left(M_T^n<x,X_T>y\right)+o\left(\!\frac{1}{\sqrt{n}}\!\right)\label{Correc1}
\\\P\left(M_T<x,X_T>y\right)&=&\P\!\left[M_T^n<x-\frac{\sigma\sqrt{T}\beta_1}{\sqrt{n}},X_T>y\right]\!\!+
\!o\left(\!\frac{1}{\sqrt{n}}\!\right)\!\!\!.\label{Correc2}
\end{eqnarray}
Therefore, we deduce from  Theorem~\ref{probcorrection}
the relations between continuous and discrete barrier options.
\begin{prop}\label{barriercorrection}
Let $X$ be a L\'evy process with generating triplet $(\gamma,\sigma^2,\nu)$ satisfying $\sigma>0$ and $\nu(\R)<\infty$, $V(H)$ be the price of a continuous option with barrier $H$, and $V^n(H)$ be the price of the corresponding discrete barrier option. We assume that the process $\left(e^{X_t-(r-\delta)t}\right)_{t\geq0}$ is a martingale. Then
\begin{eqnarray*}
V^n(H)&=&V\left(He^{\pm\frac{\sigma\sqrt{T}\beta_1}{\sqrt{n}}}\right)+o\left(\frac{1}{\sqrt{n}}\right)
\\V\left(H\right)&=&V^n\left(He^{\mp\frac{\sigma\sqrt{T}\beta_1}{\sqrt{n}}}\right)+o\left(\frac{1}{\sqrt{n}}\right),
\end{eqnarray*}
where in $\pm$ and $\mp$, the top case applies for \emph{Up} options and the bottom case applies for \emph{Down} options.
\end{prop}

\begin{rmq}\rm
Under the assumptions of Proposition~\ref{barriercorrection}, we can prove that
\begin{eqnarray*}
V(H)-V^n\left(H\right)&=&\frac{C}{\sqrt{n}}+o\left(\frac{1}{\sqrt{n}}\right).
\end{eqnarray*}
See Remark~\ref{rem:fin}.
\end{rmq}

\begin{pf} \itshape of Proposition~\ref{barriercorrection}. \upshape
For the proof, we will consider only barrier options without \emph{rebate}, since the latter can be easily deduced from the former. Theorem~\ref{probcorrection} and (\ref{Correc1}) lead obviously to the following results
\begin{eqnarray*}
&&\P\left[M_T<y+\frac{\sigma\sqrt{T}\beta_1}{\sqrt{n}},X_T\leq x\right]=\P\left[M_T^n<y,X_T\leq x\right]+o\left(\frac{1}{\sqrt{n}}\right)
\\&&\P\left[M_T\geq y+\frac{\sigma\sqrt{T}\beta_1}{\sqrt{n}},X_T\leq x\right]=\P\left[M_T^n\geq y,X_T\leq x\right]+o\left(\frac{1}{\sqrt{n}}\right)
\\&&\P\left[M_T< y+\frac{\sigma\sqrt{T}\beta_1}{\sqrt{n}},X_T<x\right]=\P\left[M_T^n< y,X_T<x\right]+o\left(\frac{1}{\sqrt{n}}\right).
\end{eqnarray*}
The price of barrier options can be written in terms of the above probabilities (as in the case of the call \emph{Up and Out} studied in the beginning of the section). Recall that in the \emph{Down} case, the infimum process $m$ of $X$ satisfies
\begin{eqnarray*}
m_t&=&\inf_{0\leq s\leq t}X_s
\\&=&-\sup_{0\leq s\leq t}(-X_s).
\end{eqnarray*}
We deduce the first result of the proposition. For the second part of the proposition, we proceed in the same way and use (\ref{Correc2}).
\end{pf}

We will test the performance of Proposition~\ref{barriercorrection} with the double exponential \emph{jump-diffusion} model (see \cite{kou02}). So, we have
\begin{equation*}
X_s=\gamma s+\sigma B_s + \sum_{i=1}^{N_s}Y_i,
\end{equation*}
where $N$ is a poisson  process with intensity $\lambda$, and $Y_1$ follows an asymmetric double exponential distribution with probability density function
\begin{equation*}\label{density}
f_Y(y)=p\eta_1 e^{-\eta_1 y}\indicatrice_{\{y\geq0\}}+q\eta_2 e^{\eta_2 y}\indicatrice_{\{y<0\}},
\end{equation*}
where $\eta_1$, $\eta_2$  are positive numbers (with $\eta_1>1$ 
to ensure integrability of the exponential), 
and the non-negative real numbers $p$ and $q$ satisfy $p+q=1$. In our numerical examples,
the values of the parameters are the following: $\sigma=0.3$, $p=0.6$, $\lambda=7$, 
$\eta_1=50$ and $\eta_2=25$. We will consider the \emph{up and out put} option with parameters 
$S_0=100$, $r=0.05$, $\delta=0$, $T=1$, $K=100$, $H=110$ and $rebate=10$.
 The continuous price,  computed by the method used in \cite{kou-wang04}, is equal to $13.240$. 
 The discrete prices are computed by Monte Carlo methods. 
 In Table~\ref{put_tab}, we study the convergence of the discrete price and 
 the corrected discrete price (using the second equality in Proposition~\ref{barriercorrection}) to the 
 continuous price.
\begin{table}[ht]
\begin{center}
\begin{tabular}{|c|c|c|c|c|}
 \hline 
 \small n&Discrete price&Relative error&Corrected discrete price&Relative error\\ \hline
 \small 5&$14.193$&$7.201\%$&$13.883$&$4.857\%$\\ \hline
 \small 6&$14.160$&$6.945\%$&$13.772$&$4.016\%$\\ \hline
 \small 7&$14.128$&$6.707\%$&$13.667$&$3.379\%$\\ \hline
 \small 8&$14.095$&$6.459\%$&$13.627$&$2.923\%$\\ \hline
 \small 9&$14.072$&$6.278\%$&$13.577$&$2.544\%$\\ \hline
 \small 10&$14.048$&$6.101\%$&$13.542$&$2.273\%$\\ \hline
 \small 15&$13.957$&$5.410\%$&$13.429$&$1.430\%$\\ \hline
 \small 25&$13.851$&$4.619\%$&$13.358$&$0.896\%$\\ \hline
\end{tabular}
\caption{Performance of the continuity correction in double exponential jump-diffusion model.}
\label{put_tab}
\end{center}
\end{table}

As expected, the discrete price converges slowly, while the corrected price converges 
rapidly to the continuous price.
The reverse problem is studied in Table~\ref{kousumbar}. 
We approximate the discrete barrier price by the corrected continuous price according to our
correction formula (see the first result of Proposition~\ref{barriercorrection}). 
In the last column we give the relative error made by approximating the discrete price by the corrected 
continuous price. The latter clearly is a good approximation of the discrete price, 
compared to the continuous price.

\begin{table}[ht]
\begin{center}
\begin{tabular}{|c|c|c|c|}
 \hline 
 \small n&Discrete price&Corrected continuous price& Relative Error\\ \hline
  \small 5 &$14.193$&$13.964$&$1.613\%$\\ \hline 
  \small 10&$14.048$&$13.894$&$1.096\%$\\ \hline 
  \small 15&$13.957$&$13.829$&$0.917\%$\\ \hline 
  \small 20&$13.896$&$13.780$&$0.834\%$\\ \hline 
  \small 25&$13.851$&$13.742$&$0.787\%$\\ \hline 
  \small 30&$13.816$&$13.711$&$0.760\%$\\ \hline 
  \small 35&$13.789$&$13.685$&$0.754\%$\\ \hline 
  \small 40&$13.764$&$13.664$&$0.726\%$\\ \hline 
  \small 45&$13.743$&$13.645$&$0.713\%$\\ \hline 
  \small 50&$13.726$&$13.629$&$0.707\%$\\ \hline 
\end{tabular}
\caption{Performance of the continuity correction in double exponential jump-diffusion model.}
\label{kousumbar}
\end{center}
\end{table}
\section{Estimates for the Poisson process}\label{sec:Poisson}
%
%
In this section, we give some estimates for the jump times of a Poisson process.
These estimates will be used to derive  domination conditions in order to justify
the convergence of some expectations.
\begin{prop}\label{prop-temps}
Let $(N_t)_{t\geq 0}$ be a homogeneous Poisson process, with jump times $(T_l)_{l\geq 1}$. For $t>0$ fixed and for any integer $l\geq 1$, we have, for $i=1,\ldots, l$,
\[
\E\left(\frac{1}{\sqrt{T_{i}-T_{i-1}}}\;|\; N_t=l\right)\leq \frac{2l}{\sqrt{t}}
\]
and
\[
\E\left(\frac{1}{\sqrt{t-T_l}}\;|\; N_t=l\right)\leq \frac{2l}{\sqrt{t}}.
\]
\end{prop}

\begin{pf}.
Using the conditional distribution of the jump times $T_1$,\ldots, $T_l$, given $\{N_t=l\}$, we have
\begin{eqnarray*}
\E\left(\frac{1}{\sqrt{T_{i}-T_{i-1}}}\;|\;N_t=l\right)& = & \int_{\{0<t_1<\ldots<t_l<t\}}\frac{1}{\sqrt{t_i-t_{i-1}}}\frac{l!}{t^l}dt_1\ldots dt_l \\
& = & \int_{\left\{u_1>0,\ldots,u_l>0,\sum_{j=1}^lu_j<t\right\}}
                           \frac{1}{\sqrt{u_{i}}}\frac{l!}{t^l}du_1\ldots du_l \\
&\leq &\frac{2l}{\sqrt{t}}\int_{\left\{u_1>0,\ldots,u_l>0,\sum_{j\neq i} u_j<t\right\}}
                    \frac{(l-1)!}{t^{l-1}}du_1\ldots\mbox{\sl\strike{du}}_i \ldots du_l\\
&= &\frac{2l}{\sqrt{t}},
\end{eqnarray*}
where we have used $\int_0^t\frac{du_i}{\sqrt{u_i}}=2/\sqrt{t}$ and, in the last integral, $u_i$ is omitted. The proof of the second inequality is similar.
\end{pf}

\begin{prop}\label{prop-temps2}
Let $(N_t)_{t\geq 0}$ be a homogeneous Poisson process, with jump times $(T_l)_{l\geq 1}$. For $t>0$ fixed and for any integer $l\geq 1$, we have, for $i=1,\ldots, l$ and for any $\alpha>0$,
\[
\P\left(T_i-T_{i-1}\leq \alpha t\;|\;N_t=l\right)\leq l\alpha
\]
and
\[
\P\left(t-T_l\leq \alpha t\;|\;N_t=l\right)\leq l\alpha.
\]
\end{prop}

\begin{pf}.
We can assume that $\alpha<1$ and write, for $i=1,\ldots,l$, using the conditional distribution of jump times given
$\{N_t=l\}$,
\begin{eqnarray*}
\P\left(T_{i}-T_{i-1}\leq \alpha t\;|\;N_t=l\right)& = & 
           \int_{\left\{0<t_1<\ldots<t_l<t\right\}}
           \indicatrice_{\{t_i-t_{i-1}\leq \alpha t\}}\frac{l!}{t^l}dt_1\ldots dt_l \\
& = & \int_{\left\{u_1>0,\ldots,u_l>0,\sum_{j=1}^lu_j<t\right\}}
      \indicatrice_{\{u_{i}\leq \alpha t\}}\frac{l!}{t^l}du_1\ldots du_l \\
&\leq &l\alpha\int_{\left\{u_1>0,\ldots,u_l>0,\sum_{j\neq i} u_j<t\right\}}
            \frac{(l-1)!}{t^{l-1}}du_1\ldots\mbox{\sl\strike{du}}_i \ldots du_l\\
&= &l\alpha,
\end{eqnarray*}
where, in the last integral, the variable $u_i$ is omitted. The proof of the second inequality is similar.
\end{pf}

\section{Conditioning with respect to the jump times}\label{sec:jump}
For the proof of Theorem~\ref{probcorrection}, we will first condition
with respect to the jump times of the Poisson process.
Fix $x>0$ and $t>0$. We have
\begin{eqnarray*}
\E\left(g(X_t)\indicatrice_{\{M_t\geq x\}}\right)-
   \E \left(g(X_t)\indicatrice_{\{M^n_t\geq x\}}\right) & = &
              \E \left(g(X_t)\indicatrice_{\{M_t\geq x>M^n_t\}}\right) \\
 & = & \sum_{l=0}^\infty\E \left(g(X_t)
          \indicatrice_{\{M_t\geq x>M^n_t\}}\;|\; N_t=l\right)\P(N_t=l).
\end{eqnarray*}
Conditionally on $\{N_t=0\}$, we have $X_s=\gamma s +\sigma B_s$, for $s\in[0,t]$,
and
\[
\E\left(g(X_t)\indicatrice_{\{M_t\geq x>M^n_t\}}\;|\; N_t=0\right)=
        \E\left(g(X_t)\indicatrice_{\{M^0\geq x>M^{0,n}\}}\;|\; N_t=0\right),
\]
with $M^0=\sup_{0\leq s\leq t}(\gamma s+\sigma B_s)$, $M^{0,n}=\max_{k=0,\ldots,n}X_{kt/n}$.

Conditionally on $\{N_t=l\}$ and $\{T_1=t_1,\ldots,T_l=t_l\}$, with $0<t_1<\ldots<t_l<t$, we set $t_{l+1}=t$ and, for $j=0,\ldots,l$,
\[
M^j=\sup_{ s\in[t_j,t_{j+1})}X_s,\quad M^{j,n}=\max_{k, kt/n\in [t_j,t_{j+1})}X_{kt/n},
\]
with, by convention $M^{j,n}=-\infty$ if there is no integer $k$ such that $kt/n\in [t_j,t_{j+1})$. 
In the sequel,  we denote by $\theta$ the vector $(t_1,\ldots,t_l)$ and by $\E_{l,\theta}$ 
the conditional expectation given $\{N_t=l,T_1=t_1,\ldots,T_l=t_l\}$.
Conditionally on the values of $X$ at times $t_j$, the random variables $M^j$ 
are independent and have a density. So they are almost surely pairwise distinct and we have
\begin{eqnarray*}
\E_{l,\theta}\left(g(X_t)\indicatrice_{\{M_t\geq x>M^n_t\}}\right) & = & 
              \sum_{j=0}^l\E_{l,\theta}\left(g(X_t)
              \indicatrice_{\{M_t\geq x>M^n_t , M^j>\max_{i\neq j}M^i\}}\right)\\
 & = &  \sum_{j=0}^l\E_{l,\theta}\left(g(X_t)
              \indicatrice_{\{M^j\geq x>M^n_t , M^j>\max_{i\neq j}M^i\}}\right).
\end{eqnarray*}
Hence
\[
\E_{l,\theta}\left(g(X_t)\indicatrice_{\{M_t\geq x>M^n_t\}}\right) =\sum_{j=0}^l \left(\alpha_l^{j,n}(\theta)-\beta_l^{j,n}(\theta)\right),
\]
with
\[
\alpha_l^{j,n}(\theta)=
     \E_{l,\theta}\left(g(X_t)
              \indicatrice_{\{M^j\geq x>M^{j,n} , M^j>\max_{i\neq j}M^i\}}\right)
\]
and
\[
\beta_l^{j,n}(\theta)=
          \E_{l,\theta}\left(g(X_t)
              \indicatrice_{\{M^j\geq x>M^{j,n} , M^j>\max_{i\neq j}M^i, \max_{i\neq j}M^{i,n}\geq x\}}\right).
 \]
Integrating with respect to the jump times, we get
\[
\E\left(g(X_t)
  \indicatrice_{\{M_t\geq x>M^n_t\}}\right)=\E\left(\alpha^n_{N_t}(T_1,\ldots,T_{N_t})-\beta^n_{N_t}(T_1,\ldots,T_{N_t})\right),
\]
where for $l\in\N$,
\[
\alpha^n_{l}(t_1,\ldots,t_l)=\indicatrice_{\{l\geq 1\}}\sum_{j=0}^l \alpha^{j,n}_{l}(t_1,\ldots,t_l)
\]
and
\[
\beta^n_{l}(t_1,\ldots,t_l)=\indicatrice_{\{l\geq 1\}}\sum_{j=0}^l \beta^{j,n}_{l}(t_1,\ldots,t_l).
\]
With these notations, we can state the following proposition.
\begin{prop}\label{prop-2.1}
We have $\lim_{n\to +\infty}\sqrt{n}\E\left(\beta^n_{N_t}(T_1,\ldots,T_{N_t})\right)=0$.
\end{prop}

For the proof of  this proposition, we will use the following reformulation
 of the Asmussen-Glynn-Pitman Theorem. It can be deduced from a careful reading of the proof of 
 Theorem 1 in \cite{asmussen-al95} (see particularly pages 879 to 883,
and Remark 2).
\begin{thm}\label{Thm-AGP}
  Consider four real numbers $a$, $b$, $x$ and $y$, with $0\leq a<b$. Let 
  $\beta =(\beta_s)_{a\leq s\leq b}$
  be a Browian bridge from $x$ to $y$ over the time interval $[a,b]$ 
  (so that $\beta_a=x$ and $\beta_b=y$) and let $t$ be a fixed positive number.
  Denote by $M$ the supremum of $\beta$ and, for any positive integer $n$, by $M^n$
   the discrete supremum
  associated with a mesh of size $t/n$, so that
  \[
  M=\sup_{a\leq s\leq b}\beta_s\quad\mbox{and}\quad M^n=\sup_{k\in I_n}\beta_{\frac{kt}{n}},
  \mbox{ where } I_n=\left\{k\in \N\;|\; \frac{kt}{n}\in [a,b]\right\}.
  \]
  Then, as $n$ goes to infinity, the pair $\left(\sqrt{n}\left(M-M^n\right), \beta\right)$ converges in distribution to the pair $(\sqrt{t}R,\beta)$
  where $R$, defined in Theorem~\ref{probcorrection}, is independent of $\beta$.
\end{thm}

\noindent{\it Proof of Proposition~\ref{prop-2.1}.}
We have
\begin{eqnarray*}
|\beta^{j,n}_l(\theta)| & \leq& 
   ||g||_\infty\sum_{i\neq j}\P_{l,\theta}\left(M^j\geq x>M^{j,n},M^j>M^i\geq M^{i,n}\geq x\right)
 \\
  &\leq &||g||_\infty \sum_{i\neq j}\P_{l,\theta}\left(x\leq M^j< x+(M^j-M^{j,n}),
           x\leq M^i< x+(M^j-M^{j,n})\right).
\end{eqnarray*}
Conditionally on $\{N_t=l\}$, $\{(T_1,\ldots,T_l)=\theta\}$ and $\{X_{T_k}=x_k, k=1,\ldots l\}$, the processes
$(X_s-X_{t_i})_{t_i\leq s<t_{i+1}}$ and $(X_s-X_{t_j})_{t_j\leq s<t_{j+1}}$ (for $i\neq j$) are independent 
Brownian motions. The pairs of random variables $(M^j-x_j,M^j-M^{j,n})$ and $(M^i-x_i,M^i-M^{i,n})$ are
 independent and we have
\[
\E\left(\indicatrice_{\{x\leq M^i< x+M^j-M^{j,n}\}}\;|\;M^j,M^{j,n}\right)=\int_{x-x_i}^{x-x_i+M^j-M^{j,n}}f_i(u)du,
\]
where $f_i$ is the probability density function of the random variable \\$\sup_{0\leq s\leq t_{i+1}-t_i}
\left(\gamma s+\sigma B_s\right)$. We know (see for example Lemma 2.22 of \cite{dia}) that the function 
$f_i$ is bounded by $C/\sqrt{t_{i+1}-t_i}$, where the constant $C$ depends only on $\gamma$, 
$\sigma$ and $t$. We deduce that
\begin{eqnarray*}
&&\P_{l,\theta}\left(x\leq M^j< x+(M^j-M^{j,n}),x\leq M^i< x+(M^j-M^{j,n})\right) 
\\&& \leq \frac{C}{\sqrt{t_{i+1}-t_i}}\E_{l,\theta}\left((M^j-M^{j,n})\indicatrice_{\{x\leq M^j< x+(M^j-M^{j,n})\}}\right)\!,
\end{eqnarray*}
Note that by 
Theorem~\ref{Thm-AGP} and the fact that $M^j$ has a continuous distribution, the sequence
$\left(\sqrt{n}(M^j-M^{j,n})\indicatrice_{\{x\leq M^j< x+(M^j-M^{j,n})\}}\right)_{n\in\N}$ converges to $0$
in probability and, since $\left(\sqrt{n}(M^j-M^{j,n})\right)_{n\in\N}$ is uniformly
integrable (see \cite{asmussen-al95}, Lemma 6),
we have
\begin{eqnarray*}
\lim_{n\to +\infty}\sqrt{n}\E_{l,\theta}\left((M^j-M^{j,n})\indicatrice_{\{x\leq M^j< x+(M^j-M^{j,n})\}}\right)=0.
\end{eqnarray*}
On the other hand, we have
\begin{eqnarray}
\sum_{j=0}^l \beta^{j,n}_l(\theta)&\leq &\sum_{j=0}^l \sum_{i\neq j}\frac{C}{\sqrt{t_{i+1}-t_i}}\sqrt{\frac{t}{n}} \\
 & \leq& Cl\sqrt{\frac{t}{n}}\sum_{i=0}^l\frac{1}{\sqrt{t_{i+1}-t_i}}.
\end{eqnarray}
We deduce that the sequence of random variables$\sqrt{n}\beta^n_{N_t}(T_1,\ldots,T_{N_t})$ is dominated by an integrable random variable by Proposition~\ref{prop-temps}. This concludes the proof.
\hfill\mbox{$\diamond$}\medskip

\section{Conditioning with respect to the positions at jump times
and representation using the Bessel process}
\label{section-rep}

It follows from the discussion of Section~\ref{sec:jump} that
\begin{eqnarray*}
\E\left(g(X_t)
  \indicatrice_{\{M_t\geq x>M^n_t\}}\right)&=&
  \E_0\left(g(X_t)
  \indicatrice_{\{M^0\geq x>M^{0,n}_t\}}\right)\P(N_t=0)
  \\
  &&+\; \E\left(\alpha^n_{N_t}(T_1,\ldots,T_{N_t})\right)+o(1/\sqrt{n}),
\end{eqnarray*}
where $\E_0=\E\left(.\;|\;N_t=0\right)$ and $\alpha^n_{l}(t_1,\ldots,t_l)=\sum_{j=0}^l \alpha^{j,n}_{l}(t_1,\ldots,t_l)$.
We have 
\[
\alpha^{j,n}_l(t_1,\ldots,t_l)=\E\left(g(X_t)
  \indicatrice_{\{ M^j\geq x>M^{j,n},M^j>\max_{i\neq j}M^i\}}\;|\;N_t=l,T_1=t_1,\ldots,T_l=t_l\right).
\]

For $j=0,\ldots,l$, we set
\[
\hat{\beta}^j_u=\frac{X_{t_j+u(t_{j+1}-t_j)}}{\sigma{\sqrt{t_{j+1}-t_j}}},\mbox{ for } u\in[0,1),\quad
     \mbox{ and } \hat{\beta}^j_1=\frac{X_{t^-_{j+1}}}{\sigma{\sqrt{t_{j+1}-t_j}}}.
\]
We have, putting $\sigma_j=\sigma{\sqrt{t_{j+1}-t_j}}$,
\[
M^j=\sigma_j\hat{M}^j\quad\mbox{and}\quad M^{j,n}=\sigma_j\hat{M}^{j,n},
\]
where
\[
\hat{M}^j=\sup_{u\in[0,1]}\hat{\beta}^j_u \quad\mbox{and}\quad \hat{M}^{j,n}=\sup_{k\in I^j_n}\hat{\beta}^j_{\frac{\lambda_j k}{n}-\hat{t}_j},
\]
with the notations $\lambda_j=t/(t_{j+1}-t_j)$, $\hat{t}_j=t_j/({t_{j+1}-t_j})$ and 
\[
I^j_n=\{k\in \N\;|\; t_j\leq kt/n<t_{j+1}\} \mbox{ for }j=0,\ldots,l-1
\]
and $I^l_n=\{k\in \N\;|\; t_j\leq kt/n\leq t_{l+1}=t\}$. Here again we use
the convention $\hat{M}^{j,n}=-\infty$ if $I^j_n=\emptyset$.

For the computation of  $\alpha^{j,n}_l(t_1,\ldots,t_l)$, we will further condition with respect to 
$\{X_{T_1}=x_1,\ldots,X_{T_l}=x_l\}$, where $x_1$, \ldots, $x_l$ are arbitrary real numbers. 
So, we introduce the notations
\[
\theta=(t_1,\ldots,t_l), \quad
\xi=(x_1,\ldots,x_l),
\]
and
\[
\P_{l,\theta,\xi}=\P\left(\cdot \;|\;N_t=l,T_k=t_k,X_{T_k}=x_k,k=1,\ldots,l \right).
\]
The expectation under $\P_{l,\theta,\xi}$ will be denoted by $\E_{l,\theta,\xi}$.
Note that, under $\P_{l,\theta,\xi}$, the processes $\hat{\beta}^j$ are independent Brownian motions
on the interval $[0,1]$, with initial values $\hat{\beta}^j_0=\hat{x}_j$, with $\hat{x}_j=x_j/\sigma_j$,
$j=0, \ldots,l$.

With these notations, we can state the following lemma.
\begin{lemme}\label{lem-*}We have
\[
\E_{l,\theta,\xi}\left(g(X_t)
  \indicatrice_{\{ M^j\geq x>M^{j,n},M^j>\max_{i\neq j}M^i\}}\right)
=\E_{l,\theta,\xi}\left(
  \indicatrice_{\{ M^j\geq x>M^{j,n}\}}
  \alpha^j_{l,\theta,\xi}(X_{t_{j+1}^-},M^j)\right),
\]
where
\[
\alpha^j_{l,\theta,\xi}(y,m)=
   \left\{
     \begin{array}{l}
           \E_{l,\theta,\xi}\left(g(X_t)\indicatrice_{\{\max_{i\neq j}M^i<m\}}\right)
              \mbox{ if }  j=0,\ldots, l-1\\
         g(y)\P_{l,\theta,\xi}\left(\max_{i\neq j}M^i<m\right) \mbox{ if }  j=l.
         \end{array}
        \right.
\]
\end{lemme}
\begin{pf}.
Note that, if $j<l$, under $\P_{l,\theta,\xi}$, the pair $(M^j,M^{j,n})$ on the one hand, and 
the random variables $X_t$, $M^i$ for $i\neq j$ on the other hand 
are independent, so that
\[
\E_{l,\theta,\xi}\left(g(X_t)
  \indicatrice_{\{ M^j\geq x>M^{j,n},M^j>\max_{i\neq j}M^i\}}\right)
=\E_{l,\theta,\xi}\left(
  \indicatrice_{\{ M^j\geq x>M^{j,n}\}}\alpha^j_{l,\theta,\xi}(M^j)\right),
\]
where
\[
\alpha^j_{l,\theta,\xi}(m)=\E_{l,\theta,\xi}\left(g(X_t)\indicatrice_{\{\max_{i\neq j}M^i<m\}}\right).
\]
Note that in this case the function $\alpha^j_{l,\theta,\xi}$ does not depend on $y$.
Now, if $j=l$, we have $X_t=X_{t_{j+1}^-}$ and the random variables $M^i$ for $i<l$
are independent of the pair $(X_{t_{l+1}},M^l)$, so that
\begin{equation}\nonumber
\E_{l,\theta,\xi}\left(g(X_t)
  \indicatrice_{\{ M^j\geq x>M^{j,n},M^j>\max_{i\neq j}M^i\}}\right)
=\E_{l,\theta,\xi}\left(
  \indicatrice_{\{ M^j\geq x>M^{j,n}\}}
  \alpha^j_{l,\theta,\xi}(X_{t_j+1},M^j)\right),
\end{equation}
with
\[
\alpha^j_{l,\theta,\xi}(y,m)=
         g(y)\P_{l,\theta,\xi}\left(\max_{i\neq j}M^i<m\right).
\]
\end{pf}

We will now give a representation of the expectations in Lemma~\ref{lem-*} in terms of Bessel processes.
Set $\tau^j=\sup\{u\in[0,1]\;|\; \hat{\beta}^j_u=\hat{M}^j\}$. 
Conditionally on $\tau^j=s$ and $\hat{M}^j=m$, we set  $R^j_1(u)=m-\hat{\beta}^j_{s-u}$, for $u\in [0,s]$ and $R^j_2(v)=m-\hat{\beta}^j_{s+v}$, for $v\in[0,1-s]$. We know that, conditionally on $\{\tau^j=s,\hat{M}^j=m,\hat{\beta}^j_1=y\}$, the processes $R^j_1$ et $R^j_2$ are independent Bessel bridges of dimension 3 (cf. \cite{asmussen-al95},
Proposition 2). We can write, conditionally on $\{\tau^j=s,\hat{M}^j=m\}$, 
\[
\hat{M}^j-\hat{M}^{j,n}=\min_{k\in I^-_n}R^j_1(s+\hat{t}_j-\lambda_j(k/n))\wedge
  \min_{k\in I^+_n}  R^j_2(\lambda_j(k/n)-\hat{t}_j-s),
\]
with
\[
I^-_n=\{k\in\N\;|\;0\leq s+\hat{t}_j-\lambda_j(k/n)\leq s\},\quad I^+_n=\{k\in\N\;|\; s\leq \lambda_j(k/n)-\hat{t}_j\leq 1\}.
\]
Hence
\begin{equation}\label{Mj-Mjn}
\hat{M}^j-\hat{M}^{j,n}=
\min_{0\leq k\leq N^1_n}R^j_1(d^j_n(s)+\lambda_j(k/n))\wedge
  \min_{1\leq k\leq N^2_n}  R^j_2(\lambda_j(k/n)-d^j_n(s)),
\end{equation}
with $d^j_n(s)=\hat{t}_j+s-\frac{\lambda_j}{n}\left[\frac{n(\hat{t_j}+s)}{\lambda_j}\right]$ (where $[x]$
 is the greatest integer in $x$;
note that $0\leq d^j_n(s)\leq \lambda_j/n$) and 
\begin{eqnarray*}
N^1_n& = & \max\{k\in\N\;|\; d^j_n(s)+\lambda_j(k/n)\leq s\}, \\
N^2_n & = & \max\{k\in\N\;|\; -d^j_n(s)+\lambda_j(k/n)\leq 1-s\}.
\end{eqnarray*}
Note that $N^1_n$ is a well defined non-negative  integer if $I^-_n$ is not empty and $N^2_n$ 
is a well defined positive integer if $I^+_n$ is not empty. When one of the two sets is empty and not the 
other, the minimum in \eqref{Mj-Mjn} is considered on the non-empty set. Note that,
if $\lambda_j/n<1$ (or, equivalently, $t_{j+1}-t_j>t/n$), at least one of the two sets is non-empty.
In the following 
proposition, we will also use the notation
\[
\gamma_j=\gamma\sqrt{t_{j+1}-t_j}/\sigma.
\]
It should be emphasized that, in the next statement,
 there is no conditioning on the terminal values of the Bessel 
 processes, in contrast with the statement of Proposition 2 of 
\cite{asmussen-al95}.
\begin{prop}\label{prop-rep-Bessel} Assume $\lambda_j/n<1$.
For any bounded Borel measurable function $F:\R^3\to \R$, we have 
\[
\E_{l,\theta,\xi}
         F\left(X_{t_{j+1}^-},M^j,M^j-M^{j,n}\right)
 = \int_0^1 ds\E\left(
           \hat{F}\left(R_1(s),R_2(1-s),R^{j,n}_s\right)\right),
\]
where
\[
\hat{F}(r_1,r_2,\rho)
=\frac{1}{2}\frac{e^{\gamma_j(r_1-r_2)-(\gamma_j^2/2)}}{r_1r_2}
 F(x_j+\sigma_j(r_1-r_2),x_j+\sigma_jr_1,\sigma_j \rho),
\]
 $R_1$ et $R_2$ are independent three-dimensional Bessel processes, starting from $0$, and
\[
R^{j,n}_s=\displaystyle\min_{-N^2_n\leq k\leq N^1_n}\check{R}(d^j_n(s)+\lambda_j(k/n)),
\]
with $\check{R}(u)=R_1(u)$ for $u\geq 0$ and $\check{R}(u)=R_2(-u)$ for $u<0$.
\end{prop}
\begin{pf}.
Note that $F\left(X_{t_{j+1}^-},M^j,M^j-M^{j,n}\right)=
  F\left(\sigma_j\hat{\beta}_1,\sigma_j\hat{M}^j,\sigma_j(\hat{M}^j-\hat{M}^{j,n})\right)$.
In view of the discussion before the statement of 
Proposition~\ref{prop-rep-Bessel}, we observe that the 
conditional distribution of $\hat{M}^j-\hat{M}^{j,n}$
given $\{\tau^j=s, \hat{M}^j= m, \hat{\beta}^j_1=y\}$ is the same
as the conditional distribution of $R^{j,n}_s$ given
$\{R_1(s)=m-\hat{x}_j, R_2(1-s)=m-y\}$, so that
\begin{eqnarray*}
\lefteqn{
\E_{l,\theta,\xi}\left(
F\left(X_{t_{j+1}^-},M^j,M^j-M^{j,n}\right)
  \;|\; \tau^j=s, \hat{M}^j= m, \hat{\beta}^j_1=y\right)=}\\
&&\E\left(
    F\left(\sigma_jy,\sigma_jm,\sigma_jR^{j,n}_s\right)
    \;|\; R_1(s)=m-\hat{x}_j, R_2(1-s)=m-y
    \right)\\
    &=&
  \psi^j_s(m-\hat{x}_j,m-y),
\end{eqnarray*}
with
\begin{eqnarray*}
 \psi^j_s(r_1,r_2)= \E\left(
 G\left(r_1,r_2,R^{j,n}_s\right)
    \;|\; R_1(s)=r_1, R_2(1-s)=r_2
    \right),
\end{eqnarray*}
and
\[
 G(r_1,r_2,\rho)=F\left(x_j+\sigma_j(r_1-r_2), x_j+\sigma_j r_1, \sigma_j\rho\right).
 \]

Recall that, under probability $\P_{l,\theta,\xi}$, the process 
$(\hat{\beta}^j_u)_{u\in[0,1]}$ is a Brownian motion, starting  
from $\hat{x}_j= x_j/\sigma_j$, with drift $\gamma_j$ and with unit 
variance coefficient. It follows that the conditional probability density function of 
the pair $(\tau^j,\hat{M}^j)$ given $\hat{\beta}^j_1=y$ can be written
\[
\P(\tau^j\in ds,\hat{M}^j\in dm\;|\; \hat{\beta}^j_1=y)=\frac{u(s,m-\hat{x}_j)u(1-s,m-y)}{n(y-\hat{x}_j)}dsdm,
\]
where $n$ is the probability density function of the standard normal distribution and
\[
u(s,m)=\frac{1}{\sqrt{\pi}}\frac{m}{s^{3/2}}e^{-m^2/(2s)}, \quad m>0.
\]
The above expression of the conditional distribution of 
$(\tau^j,\hat{M}^j)$ follows from Proposition $8.15$
 in Chapter II of \cite{karatzas-shreve}. 
 
 We now have 
 \begin{eqnarray*}
 \lefteqn{
\E_{l,\theta,\xi} F\left(X_{t_{j+1}^-},M^j,M^j-M^{j,n}\right)=}\\
&&
\int
  \P(\tau^j\in ds,\hat{M}^j\in dm,\hat{\beta}^j_1\in dy)
     \psi^j_s(m-\hat{x}_j,m-y)\\
     &=&
  \int \P(\hat{\beta}^j_1\in dy)\int_0^1ds
  \int_{\hat{x}_j\vee y}^\infty \frac{dm}{n(y-\hat{x}_j)}
   u(s,m-\hat{x}_j)u(1-s,m-y)  \psi^j_s(m-\hat{x}_j,m-y).
\end{eqnarray*}            
Since $\P(\hat{\beta}^j_1\in dy)=n\left(y-\hat{x}_j-\gamma_j\right)dy$, we can write, 
with the substitution $r_1=m-\hat{x}_j$, $r_2=m-y$ in the integral with respect to $y$ and $m$,
\begin{eqnarray}
\lefteqn{\E_{l,\theta,\xi}F\left(X_{t_{j+1}^-},M^j,M^j-M^{j,n}\right)=}
  \nonumber\\
 &  &\int_0^1ds
  \int_{0}^{\infty} dr_1
  \int_{0}^\infty dr_2
  e^{ \gamma_j(r_1-r_2)-(\gamma_j^2/2)}
   u(s,r_1)u(1-s,r_2)  \psi^j_s(r_1,r_2).\label{**}
   \end{eqnarray}
 Recall that the transition density of the three-dimensional Bessel process is 
 given by
\[
\tilde{q}_t(x,y)=\frac{1}{x}q_t(x,y)y,\quad x,y>0, t>0,
\]
where $q_t(x,y)$ is the transition density of Brownian motion (on $[0,+\infty)$) killed when it hits $0$, which can be written
\[
q_t(x,y)=g_t(x-y)-g_t(x+y),
\]
where $g_t$ is the  density of the normal distribution with mean $0$ and variance $t$. For these properties of the Bessel process, see \cite{revuz-yor} (Chapter VI, Section $3$). For $x=0$, we have
\[
\tilde{q}_t(0,y)=\sqrt{\frac{2}{\pi}}\frac{y^2}{t^{3/2}}e^{-y^2/(2t)}, \quad y>0,\quad t>0.
\]
Note that, for any $m>0$ and for any $s>0$,
\[
u(s,m)=\frac{1}{m\sqrt{2}}\tilde{q}_s(0,m).
\]
Hence, \eqref{**} can be written as follows
\begin{eqnarray*}
\lefteqn{\E_{l,\theta,\xi}F\left(X_{t_{j+1}^-},M^j,M^j-M^{j,n}\right)=}\\
  &&\frac{1}{2}\int_0^1ds\int_0^\infty \int_0^\infty 
     \P(R^j_1(s)\in dr_1,R^j_2(1-s)\in dr_2)
     \frac{e^{ \gamma_j(r_1-r_2)-(\gamma_j^2/2)}}{r_1r_2}
      \psi^j_s(r_1,r_2)\\
      &=&
       \int_0^1ds\E\left(L^j_s
         \E\left(G(R_1(s),R_2(1-s),R^{j,n}_s)\;|\;  R_1(s),R_2(1-s)\right)      \right),
\end{eqnarray*}
where 
\begin{equation}\label{eq-L}
L^j_s=\frac{e^{ \gamma_j(R_1(s)-R_2(1-s))-(\gamma_j^2/2)}}{2R_1(s)R_2(1-s)}.
\end{equation}
The proposition then follows from the equality
\[
\hat{F}(r_1,r_2,\rho)=
\frac{1}{2}\frac{e^{\gamma_j(r_1-r_2)-(\gamma_j^2/2)}}{r_1r_2}
 G(r_1,r_2, \rho).
 \]
\end{pf}

\section{Transition density of the Bessel process}\label{sec:Bessel}
In this section, we give some estimates for the transition density $(\tilde{q}_t(x,y))_{t>0}$, $x,y\geq 0$ of the three-dimensional Bessel process. As noted previously, we have
\[
\tilde{q}_t(x,y)=\frac{1}{x}q_t(x,y)y,\quad x,y>0, t>0,
\]
where $q_t(x,y)$ is the density of the Brownian motion killed at $0$, which can be written
\[
q_t(x,y)=g_t(x-y)-g_t(x+y),
\]
where $g_t$ is the density of the normal distribution with mean $0$ and variance $t$. 
For $x=0$, we have
\[
\tilde{q}_t(0,y)=\sqrt{\frac{2}{\pi}}\frac{y^2}{t^{3/2}}e^{-y^2/(2t)}, \quad y>0,t>0.
\]
We set, for $r> 0$, $m> 0$,
\[
\bar{q}_t(r,m)=\frac{\tilde{q}_t(r,m)}{m}=\frac{1}{r}q_t(r,m).
\]
Note that
\begin{eqnarray*}
\bar{q}_t(r,m) & = & \frac{1}{r}\left(g_t(r-m)-g_t(r+m)\right)\\
  & = &   -\int_{-1}^{+1}g'_t(m+r\xi)d\xi.
\end{eqnarray*}
This last expression allows to extend by continuity the definition of 
$\bar{q}_t(r,m)$ for $r=0$ or $m=0$, so that
\[
\bar{q}_t(0,m)=-\int_{-1}^{+1}g'_t(m)d\xi=-2g'_t(m)=\sqrt{\frac{2}{\pi}}\frac{m}{t^{3/2}}e^{-m^2/(2t)}.
\]
Notice that, for any $r\geq 0$, $\bar{q}_t(r,0)=0$.
\begin{prop}\label{prop-barq}
We have the following estimates, for any $\gamma\in\R$.
\begin{enumerate}
  \item For any $s$, $r$, $m>0$, 
  \[
  \bar{q}_s(r,m)e^{\gamma m-(\gamma^2 s/2)}\leq \frac{1}{s}e^{\gamma_+ r}
      \left(C_1+C_2\gamma_+\sqrt{s}\right),
  \]
  with $C_1=\sqrt{\frac{2}{\pi e}}$ and $C_2=\sqrt{2/\pi}$.
  \item For any $r,m>0$,
  \[
 \int_0^1 ds e^{\gamma m-(\gamma^2 s/2)}\bar{q}_s(r,m)\leq 2^{3/2}e^{r\gamma_+}.
  \]
  \item For any $s$, $r>0$,
  \[
  \int_0^{+\infty} dm e^{\gamma m-(\gamma^2 s/2)}\bar{q}_s(r,m)\leq 
      \frac{4\sqrt{2}}{\sqrt{\pi s}}e^{r\gamma_+}.
  \]
\end{enumerate}
\end{prop}

\begin{pf}.
Note that we can assume that $\gamma\geq 0$ because, for $\gamma<0$, $e^{\gamma m-(\gamma^2 s/2)}\leq 1$.
We have, using the equalities $g_t(x)=g_1(x/\sqrt{t})/\sqrt{t}$ and $g'_1(x)=-xg_1(x)$
\begin{eqnarray*}
\bar{q}_s(r,m) & = &-\int_{-1}^{+1}g'_s(m+r\xi)d\xi\\
  & = &   \frac{1}{s}\int_{-1}^{+1}\frac{m+r\xi}{\sqrt{s}}g_1\left(\frac{m+r\xi}{\sqrt{s}}\right)d\xi.
\end{eqnarray*}
Note that
\[
e^{\gamma m-(\gamma^2 s/2)}g_1\left(\frac{m+r\xi}{\sqrt{s}}\right) =
e^{-\gamma r \xi}g_1\left(\frac{m+r\xi-\gamma s}{\sqrt{s}}\right).
\]
Hence
\begin{eqnarray}
e^{\gamma m-(\gamma^2 s/2)}\bar{q}_s(r,m)&=&\frac{1}{s}\int_{-1}^{+1}e^{-\gamma r \xi}
\frac{m+r\xi}{\sqrt{s}}g_1\left(\frac{m+r\xi-\gamma s}{\sqrt{s}}\right)d\xi \label{1}\\
   &\leq &\frac{e^{r\gamma}}{s}\int_{-1}^{+1}
\left|\frac{m+r\xi-\gamma s}{\sqrt{s}}+\gamma \sqrt{s}\right|g_1\left(\frac{m+r\xi-\gamma s}{\sqrt{s}}\right)d\xi\nonumber \\
   &\leq &2\frac{e^{r\gamma}}{s}\left(\sup_{x>0}xg_1(x)+\gamma\sqrt{s}\frac{1}{\sqrt{2\pi}}\right),\nonumber
\end{eqnarray}
which gives the first inequality. For the second and third inequality, we start from (\ref{1}) and notice that
\begin{eqnarray*}
(m+r\xi-\gamma s)^2& = & (m+r\xi)^2+\gamma^2 s^2-2\gamma s(m+r\xi) \\
  & \geq & (m+r\xi)^2+\gamma^2 s^2-\left(2\gamma^2s^2+ \frac{(m+r\xi)^2}{2}\right)\\
  &=&\frac{(m+r\xi)^2}{2}-\gamma^2 s^2.
\end{eqnarray*}
We deduce
\[
g_1\left(\frac{m+r\xi-\gamma s}{\sqrt{s}}\right)\leq e^{\gamma^2 s/2}g_1\left(\frac{m+r\xi}{\sqrt{2s}}\right).
\]
Hence
\begin{equation}
\label{2}
e^{\gamma m-(\gamma^2 s/2)}\bar{q}_s(r,m)\leq 
e^{\gamma r }\int_{-1}^{+1}
\frac{|m+r\xi|}{s^{3/2}}g_1\left(\frac{m+r\xi}{\sqrt{2s}}\right)d\xi,
\end{equation}
and, integrating with respect to $s$,
\begin{eqnarray*}
\int_0^1 ds e^{\gamma m-(\gamma^2 s/2)}\bar{q}_s(r,m)&\leq&e^{\gamma r }
\int_{-1}^{+1}\left(\int_{0}^{1}
\frac{|m+r\xi|}{s^{3/2}}g_1\left(\frac{m+r\xi}{\sqrt{2s}}\right)ds\right)d\xi\\
   &= &2^{3/2}e^{r\gamma}\int_{-1}^{+1}
\left(\int_{|m+r\xi]/\sqrt{2}}^{+\infty}g_1(u)du\right)d\xi\\
    &\leq &2^{3/2}e^{r\gamma},
\end{eqnarray*}
where we have set $u=|m+r\xi|/\sqrt{2s}$. Integrating (\ref{2}) with respect to $m$, we get
\begin{eqnarray*}
\int_0^{+\infty}e^{\gamma m-(\gamma^2 s/2)}\bar{q}_s(r,m)dm&\leq& 
e^{\gamma r }\int_{-1}^{+1}\left(\int_0^{+\infty}
\frac{|m+r\xi|}{s^{3/2}}g_1\left(\frac{m+r\xi}{\sqrt{2s}}\right)dm\right)d\xi\\
& \leq &  e^{\gamma r }\int_{-1}^{+1}\left(\int_{-\infty}^{+\infty}
\frac{|m+r\xi|}{s^{3/2}}g_1\left(\frac{m+r\xi}{\sqrt{2s}}\right)dm\right)d\xi\\
&=&4\frac{e^{\gamma r }}{s^{1/2}}\int_{-\infty}^{+\infty}
|z|g_1\left(z\right)dz=4\sqrt{2}\frac{e^{\gamma r }}{\sqrt{s\pi}},
\end{eqnarray*}
where we have set $z=(m+r\xi)/\sqrt{2s}$.
\end{pf}

We complete this section with a lemma concerning the minimum of the Bessel process.
This result is a consequence of   Lemma $3$ of \cite{asmussen-al95}.
\begin{lemme}\label{lem-Bessel}
   Let $({R}(t))_{t\geq 0}$ be a three-dimensional Bessel process starting from $0$ and 
   let $t_1$, $t_2$, $y$, $m$, $b$ be positive numbers, with $t_1<t_2$. We have, using the notation 
   $R^\sharp(t_1,t_2)=\min_{u\in[t_1,t_2]}{R}(u)$,
 \[ 
   \P\left(R^\sharp(t_1,t_2)\leq b \;|\; {R}(t_1)=y,{R}(t_2)=m\right)\leq \frac{b(m+y)}{ym}.
  \] 
\end{lemme}
\begin{pf}.
We assume that $b<y\wedge m$, since if $b\geq y\wedge m$, the upper bound is larger
than  or equal to $1$. We then have, using Lemma $3$ of \cite{asmussen-al95} (and the fact that the Bessel bridge can be viewed as a Brownian bridge conditioned  to remain positive: see the 
proof of Lemma 4 of \cite{asmussen-al95}),
\[
\P\left(R^\sharp(t_1,t_2)\leq b \;|\; {R}(t_1)=y,{R}(t_2)=m\right)=
\frac{e^{2(b-y)(m-b)/T}-e^{-2ym/T}}{1-e^{-2ym/T}},
\] 
with $T=t_2-t_1$. Hence, using the convexity of the exponential function and the inequality $b(m-b+y)\leq ym$,
\begin{eqnarray*}
\P\left(R^\sharp(t_1,t_2)\leq b \;|\; {R}(t_1)=y,{R}(t_2)=m\right)&=&
\frac{e^{2((b-y)(m-b)+ym)/T}-1}{e^{2ym/T}-1}\\
&=&\frac{e^{2b(m-b+y)/T}-1}{e^{2ym/T}-1}\\
&\leq &\frac{b(m-b+y)}{ym}\leq \frac{b(m+y)}{ym}.
\end{eqnarray*}
\end{pf}

\section{Domination of the conditional probability}\label{sec:domination}
The following proposition will be used to ensure the domination of conditional expectations.
\begin{prop}\label{prop-dom}
There exists a constant $C_{t,\gamma,\sigma}$ (depending only on $t$, $\gamma$ and $\sigma$) such that, for any $l\in\N$, $\theta=(0<t_1<\ldots<t_l<t)\in\R^l$ et $\xi\in\R^l$, we have, for $j=0,\ldots, l$,
\[
\P_{l,\theta,\xi}\left(M^j\geq x>M^{j,n}\right)
\leq \indicatrice_{\{t_{j+1}-t_j\leq 8t/n\}}+
             \frac{C_{t,\gamma,\sigma}}{\sqrt{n}}\frac{1}{\sqrt{t_{j+1}-t_j}}.
\]
\end{prop}
\begin{pf}.
Note that $\P_{l,\theta,\xi}\left(M^j\geq x>M^{j,n}\right)=\E_{l,\theta,\xi}F\left(M^j,M^j- M^{j,n}\right)$,
if we define the function $F$ by $F(m,\rho)= \indicatrice_{\{x\leq m < x+\rho\}}$.
It follows from  Proposition~\ref{prop-rep-Bessel} that\begin{eqnarray*}
\P_{l,\theta,\xi}\left(M^j\geq x>M^{j,n}\right)
&= &\int_0^1 ds\E\left(L^j_s
           \indicatrice_{\left\{\tilde{x}_j\leq  R_1(s)\leq\tilde{x}_j+R^{j,n}_s\right\}}\right),
\end{eqnarray*}
where $L^j_s$ is given by \eqref{eq-L} and
\[
\tilde{x}_j=\frac{x-x_j}{\sigma_j}.
\]
We can obviously assume that $t_{j+1}-t_j>8t/n$, which, with the notations of Section~\ref{section-rep}, can be written $\lambda_j/n<1/8$ and ensures that at least one of the sets $I^-_{n}$ and $I^+_{n}$ is non-empty. So we can bound the random variable $R^{j,n}_s$ by $R^*(\lambda_j/n)$ where, for $u\in[0,1]$, we set
\[
R^*(u)=\max_{-u\leq v\leq +u}\check{R}(v).
\]
Hence
\begin{eqnarray*}
\P_{l,\theta,\xi}\left(M^j\geq x>M^{j,n}\right)
&\leq &\int_0^1 ds\E\left(L^j_s
           I_j\left(R_1(s),R^*(\lambda_j/n)\right)\right),
\end{eqnarray*}
where
\[
I_j(r,\rho)= \indicatrice_{\left\{\tilde{x}_j\leq  r\leq\tilde{x}_j+\rho\right\}}.
\]
We have
\[
\int_0^{\lambda_j/n} ds\E\left(L^j_s
           I_j\left(R_1(s),R^*(\lambda_j/n)\right)\right)
           \leq \int_0^{\lambda_j/n} ds\E\left(L^j_s\right),
\]           
and          
\[
\E\left(L^j_s\right)  = \E\left(\frac{e^{\gamma_j R_1(s)-\frac{\gamma_j^2 s}{2}}}{\sqrt{2}R_1(s)}\right)
                 \E\left(\frac{e^{\gamma_j R_2(1-s)-\frac{\gamma_j^2 (1-s)}{2}}}{\sqrt{2}R_2(1-s)}\right).
\]
By scaling, we have
\begin{eqnarray*}
\E\left(\frac{e^{\gamma_j R_1(s)-\frac{\gamma_j^2 s}{2}}}{\sqrt{2}R_1(s)}\right)&=&
      \E\left(\frac{e^{\gamma_j \sqrt{s}R_1(1)-\frac{\gamma_j^2 s}{2}}}{\sqrt{2}\sqrt{s}R_1(1)}\right)\\
      &\leq &C_{t,\gamma,\sigma}\frac{1}{\sqrt{s}}.
\end{eqnarray*}
Similarly
\[
\E\left(\frac{e^{\gamma_j R_2(1-s)-\frac{\gamma_j^2 (1-s)}{2}}}{\sqrt{2}R_2(1-s)}\right)\leq C_{t,\gamma,\sigma}\frac{1}{\sqrt{1-s}}.
\]
Therefore
\[
\int_0^{\lambda_j/n} ds\E\left(L^j_s
           \right)\leq 
           C_{t,\gamma,\sigma}\sqrt{\frac{\lambda_j}{n}}=\frac{C_{t,\gamma,\sigma}}{\sqrt{t_{j+1}-t_j}}\frac{1}{\sqrt{n}},
\]
and, by a similar argument, 
$\int_{1-\lambda_j/n}^{1} ds\E\left(L^j_s
           \right)
           \leq \frac{C_{t,\gamma,\sigma}}{\sqrt{t_{j+1}-t_j}}\frac{1}{\sqrt{n}}$.
It remains to study the integral on the interval $[\lambda_j/n,1-\lambda_j/n]$.
 Denote by $(\F_s)_{s\geq 0}$ the natural filtration of the pair $(R^1, R^2)$. 
 For $s\in [\lambda_j/n,1-\lambda_j/n]$, we have, by conditioning with respect to $\F_{\lambda_j/n}$,
\begin{eqnarray*}
&&\E\left(L^j_s \indicatrice_{\left\{R_1(s)\in[\tilde{x}_j,\tilde{x}_j+R^*(\lambda_j/n)]\right\}}\right)
\\&&=\E\left(\int_0^\infty\!\!\!dm\,\bar{q}_{s-\lambda_j/n}^j(R_1(\lambda_j/n), m)\phi^j_{1-s-\lambda_j/n}(R_2(\lambda_j/n))
                I_j(m,R^*(\lambda_j/n))\right),
\end{eqnarray*}
where 
\[
\bar{q}_{s}^j(r,m)=e^{\gamma_j m-\gamma_j^2/2}\bar{q}_s(r,m)
\quad
\mbox{and}
\quad
\phi^j_s(r)=\frac{1}{2}\int_0^\infty dm e^{-\gamma_j m}\bar{q}_s(r,m).
\]
By Proposition~\ref{prop-barq}, we have
\[
\phi^j_s(r)\leq \frac{2\sqrt{2}}{\sqrt{\pi s}}e^{r|\gamma_j|+(\gamma_j^2/2)}\leq \frac{C_{t,\gamma,\sigma}}{\sqrt{s}}e^{r|\gamma_j|}.
\]
We have
\begin{eqnarray*}
&&\int_{\lambda_j/n}^{1-\lambda_j/n}\!\!\!\!\!ds \,\bar{q}_{s-\lambda_j/n}^j(R_1(\lambda_j/n), m)\phi^j_{1-s-\lambda_j/n}(R_2(\lambda_j/n))
   \\&&=\int_0^{1-2\lambda j/n}\!\!\!\!\!ds\,\bar{q}_{s}^j(R_1(\lambda_j/n), m)\phi^j_{1-s-2\frac{\lambda_j}{n}}(R_2(\lambda_j/n)).
\end{eqnarray*}
For $s\in[0,1/2-(\lambda_j/n)]$, we have 
\[
\frac{1}{\sqrt{1-s-2\lambda_j/n}}\leq \frac{1}{\sqrt{1/2-(\lambda_j/n)}}\leq 2,
\]
because $\lambda_j/n< 1/4$. Hence
\begin{eqnarray*}
&&\int_{\lambda_j/n}^{1/2-(\lambda_j/n)}\!\!\!\!\!ds \,\bar{q}_{s-\lambda_j/n}^j(R_1(\lambda_j/n), m)\phi^j_{1-s-\lambda_j/n}(R_2(\lambda_j/n))
\\&&\leq C_{t,\gamma,\sigma}\int_0^1ds\,\bar{q}_{s}^j(R_1(\lambda_j/n), m)e^{|\gamma_j| R_2(\lambda_j/n)}\\
&&\leq C_{t,\gamma,\sigma}e^{|\gamma_j|\left(R_1(\lambda_j/n)+R_2(\lambda_j/n)\right)},
\end{eqnarray*}
where the last inequality follows from Proposition~\ref{prop-barq}. On the other hand, 
for $s\in [1/2-(\lambda_j/n),1-\lambda_j/n]$, we have, using the first inequality of 
Proposition~\ref{prop-barq},
\[
\bar{q}_{s-\lambda_j/n}^j(r, m)\leq \frac{C_{t,\gamma,\sigma}}{s-\lambda_j/n}e^{r|\gamma_j|}\leq 
   4C_{t,\gamma,\sigma}e^{r|\gamma_j|}.
\]
Here, the last inequality follows from the condition $\lambda_j/n<1/8$. We deduce
\begin{eqnarray*}
&&\int_{1/2-(\lambda_j/n)}^{1-(\lambda_j/n)}\!\!\!\!\!ds \,\bar{q}_{s-\lambda_j/n}^j(R_1(\lambda_j/n), m)\phi^j_{1-s-\lambda_j/n}(R_2(\lambda_j/n))
\\&&\leq C_{t,\gamma,\sigma}e^{|\gamma_j| R_1(\lambda_j/n)}\int_{1/2}^1ds
    \phi^j_{1-s}(R_2(\lambda_j/n))
\\&&\leq C_{t,\gamma,\sigma}e^{|\gamma_j| \left(R_1(\lambda_j/n)+R_2(\lambda_j/n)\right)}\int_{1/2}^1\frac{ds}{\sqrt{1-s}}\\
\\&&\leq C_{t,\gamma,\sigma}e^{|\gamma_j|\left(R_1(\lambda_j/n)+R_2(\lambda_j/n)\right)}.
\end{eqnarray*}
Then we have
\begin{eqnarray}
&&\int_{\lambda_j/n}^{1-\lambda_j/n}ds\E\left(L^j_s
           \indicatrice_{\left\{R_1(s)\in[\tilde{x}_j,\tilde{x}_j+R^*(\lambda_j/n)]\right\}}\right)\nonumber
           \\&&\leq C_{t,\gamma,\sigma}\E\left(e^{|\gamma_j|\left(R_1(\lambda_j/n)+R_2(\lambda_j/n)\right)}\int_0^\infty\!\!\!dm\, I_j(m,R^*(\lambda_j/n))\right)\label{19}\\
           &&\leq C_{t,\gamma,\sigma}\E\left(R^*(\lambda_j/n)e^{|\gamma_j|\left(R_1(\lambda_j/n)+R_2(\lambda_j/n)\right)}\right)\nonumber
           \\&&=C_{t,\gamma,\sigma}\sqrt{\lambda_j/n}\E\left(R^*(1)e^{|\gamma_j|\sqrt{\lambda_j/n}\left(R_1(1)+R_2(2)\right)}\right)\nonumber
           \\&&=\frac{C_{t,\gamma,\sigma}}{\sqrt{t_{j+1}-t_j}}\frac{1}{\sqrt{n}},\nonumber
\end{eqnarray}
where we have used the scaling property of the Bessel process, 
$\sqrt{\lambda_j}=\sqrt{t/(t_{j+1}-t_j)}$, and $|\gamma_j|\sqrt{\lambda_j}=|\gamma|\sqrt{t}/\sigma$.
\end{pf}
\begin{rmq}\rm\label{rmq-7.2}
It follows from the proof of the proposition that, for any $\delta>0$, we have
\[
\P_{l,\theta,\xi}\left(M^j\geq x>M^j-\delta/\sqrt{n}\right)
\leq \indicatrice_{\{t_{j+1}-t_j\leq 8t/n\}}+
             \frac{C_{t,\gamma,\sigma,\delta}}{\sqrt{n}}\frac{1}{\sqrt{t_{j+1}-t_j}}.
\]
Indeed, we have $\P_{l,\theta,\xi}\left(M^j\geq x>M^j-\delta/\sqrt{n}\right)=
   \int_0^1 ds\E\left(L^j_s
          I_j\left(R_1(s), \delta/(\sigma_j\sqrt{n})\right)\right)$,
          and we can replace $R^*(\lambda_j/n)$ with $\delta/(\sigma_j\sqrt{n})$
          in \eqref{19}.
\end{rmq}
\section{Convergence of the conditional expectation}\label{sec:convergence}
The aim of this section is to prove the following result and to deduce Theorem~\ref{probcorrection}.
\begin{thm}\label{thm-convcond}
Let $F:\R^2\to \R$ be a bounded Borel measurable function, such that $m\mapsto F(y,m)$
is continuous for all $y\in\R$. 

We have, with the notation $\E_0=\E\left(.\;|\;N_t=0\right)$,
\[
\E_0\left(F(X_t,M^0)
  \indicatrice_{\{ M^0\geq x>M^{0,n}\}}\right)=
  \E_0\left(F(X_t,M^0)
  \indicatrice_{\{ M^0\geq x>M^{0}-\sigma\beta_1\sqrt{t/n}\}}\right)+o(1/\sqrt{n}).
\]

Moreover, for 
 any positive  integer $l$, and for any $\theta=(t_1,\ldots,t_l)\in \R^l$, with $0<t_1<\ldots<t_l<t$, 
$\xi=(x_1,\ldots,x_l)\in \R^l$, 
we have,  for $j\in\{0,1,\ldots,l\}$, if $x_j\neq x$,
\begin{eqnarray*}
\E_{l,\theta,\xi}\left(F(X_{t_{j+1}^-},M^j)
  \indicatrice_{\{ M^j\geq x>M^{j,n}\}}\right)&=&
\E_{l,\theta,\xi}\left(F(X_{t_{j+1}^-},M^j)
  \indicatrice_{\{ M^j\geq x>M^{j}-\sigma\beta_1\sqrt{t/n}\}}\right)
\\
&&+o(1/\sqrt{n}),
\end{eqnarray*}
where $\beta_1$ is defined as in Theorem~\ref{probcorrection}.
\end{thm}

We will first  show how Theorem~\ref{probcorrection} can be deduced
 from Theorem~\ref{thm-convcond}.

\begin{pf} \itshape of Theorem~\ref{probcorrection}. \upshape
Observe that, with the notation $\E_0=\E(\cdot\;|\;N_t=0)$,
\begin{eqnarray*}
 \E\left(g(X_t)
  \indicatrice_{\{ M_t\geq x>M^{n}_t\}}\right)&=&\E_0\left(g(X_t)
  \indicatrice_{\{ M^0\geq x>M^{0,n}\}}\right)\P(N_t=0)\\
  &&
  +\E\left(g(X_t)
  \indicatrice_{\{ M_t\geq x>M^{n}_t\}}  \;|\; N_t\geq 1\right)\P(N_t\geq1).
\end{eqnarray*}
Using Theorem~\ref{thm-convcond}, we have
\[
\E_0\left(g(X_t)
  \indicatrice_{\{ M^0\geq x>M^{0,n}\}}\right)=
      \E_0\left(g(X_t)
  \indicatrice_{\{ M^0\geq x>M^0-\sigma\beta_1\sqrt{t/n}\}}\right)+o(1/\sqrt{n}).
\]        
On the other hand, we deduce from Proposition~\ref{prop-2.1} that 
\[
\E\left(g(X_t)
  \indicatrice_{\{ M_t\geq x>M^{n}_t\}}  \;|\; N_t\geq 1\right)
     =\E\left(\sum_{j=0}^{N_t} \alpha^{j,n}_{N_t}(T_1,\ldots,T_n)\;|\;N_t\geq 1\right)+o(1/\sqrt{n}),
\]
where, for any positive integer $l$, and  for $j=0,\ldots,l$,
\[
\alpha^{j,n}_{l}(\theta)=\E_{l,\theta}\left(
     g(X_t)
  \indicatrice_{\{ M^j\geq x>M^{j,n},M^j>\max_{i\neq j}M^i \}}\right).
\]
Note that, as a consequence of Lemma~\ref{lem-*}, we have, assuming $x_j\neq x$,
\begin{eqnarray*}
\E_{l,\theta,\xi} \alpha^{j,n}_{l}(\theta)&=&
    \E_{l,\theta,\xi}\left(
  \indicatrice_{\{ M^j\geq x>M^{j,n}\}}
  \alpha^j_{l,\theta,\xi}(X_{t_{j+1}^-},M^j)\right)\\
  &=&\E_{l,\theta,\xi}\left(
  \indicatrice_{\{ M^j\geq x>M^j-\beta_1\sigma\sqrt{t/n}\}}
  \alpha^j_{l,\theta,\xi}(X_{t_{j+1}^-},M^j)\right)+o(1/\sqrt{n})\\
  &=&
  \E_{l,\theta,\xi}\left(
  \indicatrice_{\{ M^j\geq x>M^j-\beta_1\sigma\sqrt{t/n},M^j>\max_{i\neq j}M^i \}}
  \right)+o(1/\sqrt{n}),
\end{eqnarray*}
where the second equality follows from Theorem~\ref{thm-convcond},
and the last one from the expression of $\alpha^j_{l,\theta,\xi}$ (see Lemma~\ref{lem-*}
and its proof; note that $m\mapsto\alpha^i_{l,\theta,\xi}(y,m)$ is continuous
because $\P_{l,\theta,\xi}\left(\max_{i\neq j}M_i=m\right)=0$). By taking the sum over $j=0,\ldots,l$, we deduce that, for $l$,
$\theta$ and $\xi$ fixed, we have
\begin{equation}\label{20}
\E_{l,\theta,\xi}\left(g(X_t)
  \indicatrice_{\{ M_t\geq x>M^{n}_t\}} \right)=
  \E_{l,\theta,\xi}\left(g(X_t)
  \indicatrice_{\{ M_t\geq x>M_t-\beta_1\sigma\sqrt{t/n}\}} \right)+o(1/\sqrt{n}).
\end{equation}
Observe that $\P\left(X_{T_j}=x\right)= 0$ for all jump times $T_j$ (including $T_0=0$,
since $x>0$). Therefore,
in order to get the resut for the unconditional expectations, we only need to check a domination condition.
We have
\begin{eqnarray*}
\left|\E_{l,\theta,\xi}\left(g(X_t)
  \indicatrice_{\{ M_t\geq x>M^{n}_t\}} \right)\right|&\leq &||g||_\infty
\P_{l,\theta,\xi}\left( M_t\geq x>M^{n}_t\right)\\
   &\leq &||g||_\infty\sum_{j=0}^l\P_{l,\theta,\xi}\left(M^j\geq x>M^{j,n}\right),
\end{eqnarray*}
and
\[
\P_{l,\theta,\xi}\left( M_t\geq x>M_t-\beta_1\sigma\sqrt{t/n}\right)
   \leq \sum_{j=0}^l\P_{l,\theta,\xi}\left(M^j\geq x>M^j-\beta_1\sigma\sqrt{t/n}\right).
\]
Using Proposition~\ref{prop-dom} and Remark~\ref{rmq-7.2}, we deduce that
\[
\P_{l,\theta,\xi}\left( M_t\geq x>M^{n}_t\right)
+\P_{l,\theta,\xi}\left( M_t\geq x>M_t-\beta_1\sigma\sqrt{t/n}\right)\leq  P(l,\theta),
\]
where
\[
P(l,\theta)=P_l(t_1,\ldots,t_l)=\sum_{j=0}^l \indicatrice_{\{t_{j+1}-t_j\leq 8t/n\}}+
             \frac{C_{t,\gamma,\sigma}}{\sqrt{n}}\sum_{j=0}^l \frac{1}{\sqrt{t_{j+1}-t_j}}
\]
It follows from Proposition~\ref{prop-temps2}
that
\[
\E\left(\sum_{j=0}^{N_t}\indicatrice_{\{T_{j+1}\wedge t-T_j\leq 8t/n\}}
    \right)\leq \frac{8t}{n}\E \left(N_t(N_t+1)\right),
\]
and, from Proposition~\ref{prop-temps}, that
\[
\E\left(\sum_{j=0}^{N_t}\frac{1}{\sqrt{T_{j+1}\wedge t-T_j}}\right)\leq
                                 2\frac{\E \left(N_t(N_t+1)\right)}{\sqrt{t}}.
\]
The last two inequalities are sufficient to extend the estimate \eqref{20} to unconditional expectations.
\end{pf}

For the proof of Theorem~\ref{thm-convcond}, we start from the representation given by Proposition~\ref{prop-rep-Bessel}, which reads
\[
\E_{l,\theta,\xi}\left(F(X_{t_{j+1}^-},M^j)
  \indicatrice_{\{ M^j\geq x>M^{j,n}\}}\right)
   =
\int_0^1E^{j,n}_{l,\theta,\xi}(s)ds,
\]
with
\[
E^{j,n}_{l,\theta,\xi}(s)=\E\left(L^j_s
           \alpha_{j}(R_1(s),R_2(1-s))
           \indicatrice_{\left\{\tilde{x}_j\leq  R_1(s)\leq\tilde{x}_j+R^{j,n}_s\right\}}\right),
\]
where $L^j_s=\frac{e^{ \gamma_j(R_1(s)-R_2(1-s))-(\gamma_j^2/2)}}{2R_1(s)R_2(1-s)}$,
\[
\alpha_j(r_1,r_2)=F(x_j+\sigma_j(r_1-r_2),x_j+\sigma_j r_1),
\mbox{ and } \;\tilde{x}_j=\frac{x-x_j}{\sigma_j}.
\]
Note that the function $\alpha_j$ is bounded, and $||\alpha_j||_\infty = ||F||_\infty$.
For any integer $J\geq 1$, we can write, for $n$ large enough,
\begin{equation}
\label{dec-E}
\int_0^1E^{j,n}_{l,\theta,\xi}(s)ds=\int_0^{\frac{\lambda_j (J+1)}{n}}E^{j,n}_{l,\theta,\xi}(s)ds+
\int_{1-\frac{\lambda_j J}{n}}^1E^{j,n}_{l,\theta,\xi}(s)ds+
\int_{\frac{\lambda_j (J+1)}{n}}^{1-\frac{\lambda_j J}{n}}E^{j,n}_{l,\theta,\xi}(s)ds.
\end{equation}
The first two  terms of this decomposition are controlled   via the following lemma.
\begin{lemme}\label{lem-bord}
For any integer $J\geq 1$, for any $(l,\theta,\xi)$, and for $j=0,1\ldots,l$,
we have, if $x_j\neq 0$
\[
\int_0^{\frac{\lambda_j (J+1)}{n}}|E^{j,n}_{l,\theta,\xi}(s)|ds+
\int_{1-\frac{\lambda_j J}{n}}^1|E^{j,n}_{l,\theta,\xi}(s)|ds=o(1/\sqrt{n}).
\]
\end{lemme}

\begin{pf}.
We will only consider the first integral, the argument is similar for the second term.
Note that for $n$ large enough and $s\in[0,{\lambda_j (J+1)}/{n}$], we have $1-s>\lambda_j/n$
and we can bound $R^{j,n}_s$ from above by $R_2^*(\lambda_j/n)$ 
(with $R^*_2(s)=\max_{0\leq u\leq s}R^*_2(u)$.
Then we have, using the boundedness of $F$, 
\[
\left|E^{j,n}_{l,\theta,\xi}(s)\right|\leq ||F||_\infty \E\left(\tilde{L}^j_s
     \indicatrice_{\left\{\tilde{x}_j\leq  R_1(s)\leq\tilde{x}_j+R^*_2\left(\frac{\lambda_j}{n}\right)\right\}}
     \right),
\]
where $\tilde{L}^j_s=\frac{e^{|\gamma_j|\left(R_1(s)+R_2(1-s)\right)}}{2R_1(s)R_2(1-s)}$.
Hence, with the substitution $s'=ns$,
\begin{eqnarray*}
\sqrt{n}\int_{0}^{\lambda_j(J+1)/n}\!\!\!\!\!\left|E^{j,n}_{l,\theta,\xi}(s)\right|ds
     &\leq &\sqrt{n}||F||_\infty\int_0^{\lambda_j(J+1)/n}\!\!\!\!ds
     \E\left(\tilde{L}^j_s
     \indicatrice_{\left\{\tilde{x}_j\leq  R_1(s)\leq\tilde{x}_j+R^*_2\left(\frac{\lambda_j}{n}\right)\right\}}
     \right)\\
    &=&||F||_\infty\int_0^{\lambda_j(J+1)}\!\!\!\frac{ds'}{\sqrt{n}}
        \E\left(
     \tilde{L}^j_{s'/n}\indicatrice_{\left\{\tilde{x}_j\leq  R_1(s'/n)\leq\tilde{x}_j+R^*_2\left(\frac{\lambda_j}{n}\right)\right\}}\right).
\end{eqnarray*}
By scaling, we can write, using the notation
\[
\Lambda_j(r_1,r_2,\rho)=\frac{e^{ |\gamma_j|(r_1+r_2)}}{2r_1r_2}
         \indicatrice_{\left\{\tilde{x}_j\leq  r_1\leq\tilde{x}_j+\rho\right\}},
\]
\begin{eqnarray*}
\E\left(
     \tilde{L}^j_{s/n}\indicatrice_{\left\{\tilde{x}_j\leq  R_1(s/n)\leq\tilde{x}_j+R^*_2\left(\frac{\lambda_j}{n}\right)\right\}}\right)
     &= &\E\left(
     \Lambda_j\left( R_1(s/n), R_2(1-\frac{s}{n}),R^*_2\left(\frac{\lambda_j}{n}\right)\right)\right)\\
    &=&\E\left(
     \Lambda_j\left( R_1(s)/\sqrt{n}, R_2(1-\frac{s}{n}),R^*_2\left(\frac{\lambda_j}{n}\right)\right)\right)\\
     &\leq&\sqrt{n}\E\left(
     \Lambda_j^{s,n}
     \left(R_1(s),R_2(1),
           R^*_2\left(\frac{\lambda_j}{n-s}\right)\right)\right),
\end{eqnarray*}
with 
\[
 \Lambda_j^{s,n}(r_1,r_2,\rho)=\frac{e^{ |\gamma_j|\left(\frac{r_1}{\sqrt{n}}+\sqrt{1-\frac{s}{n}}r_2\right)}}{2r_1\sqrt{1-\frac{s}{n}}r_2}
         \indicatrice_{\left\{\tilde{x}_j\leq  \frac{r_1}{\sqrt{n}}\leq\tilde{x}_j+\rho\right\}}.
\]
Hence
\begin{eqnarray*}
\sqrt{n}\int_{0}^{\lambda_j(J+1)/n}\!\!\!\!\!\left|E^{j,n}_{l,\theta,\xi}(s)\right|ds
     &\leq &||F||_\infty\int_0^{\lambda_j(J+1)}\!\!\! ds
       \E\left(
     \Lambda_j^{s,n}
     \left(R_1(s),R_2(1),
           R^*_2\left(\frac{\lambda_j}{n-s}\right)\right)\right).
\end{eqnarray*}   
Now, if $n>2\lambda_j(J+1)$, we have, for $s\in[0,\lambda_j(J+1)]$,
$n-s\geq n/2$, so that $R^*_2\left(\frac{\lambda_j}{n-s}\right)\leq R^*_2\left(\frac{2\lambda_j}{n}\right)$,
and $\sqrt{1-s/n}\geq 1/\sqrt{2}$, so that
 \begin{eqnarray*}
 \Lambda_j^{s,n}
     \left(R_1(s),R_2(1),
           R^*_2\left(\frac{\lambda_j}{n-s}\right)\right)    &\leq&
    \frac{e^{|\gamma_j|(R_1(s)+R_2(1))}}{\sqrt{2}R_1(s)R_2(1)}
          \indicatrice_{\left\{\tilde{x}_j\leq  R_1(s)/\sqrt{n}\leq\tilde{x}_j+
             R^*_2\left(\frac{2\lambda_j}{n}\right)\right\}}.   
     \end{eqnarray*}
     If $x\neq x_j$, we have $\tilde{x}_j\neq 0$, and the right hand side of the inequality goes to $0$
     almost surely  as $n\to \infty$, for all $s\in(0,1)$. Since
     $
      \int_0^1ds \E\left( \frac{e^{|\gamma_j|(R_1(s)+R_2(1))}}{\sqrt{2}R_1(s)R_2(1)}\right)<\infty,
      $
      we conclude that 
      \[
      \lim_{n\to \infty}\left(\sqrt{n}\int_{0}^{\lambda_j(J+1)/n}\!\!\!\!\!E^{j,n}_{l,\theta,\xi}(s)ds \right)=0.
      \]
\end{pf}

We will now examine the case $x<x_j$. 
\begin{lemme}\label{lem-8.3-}
If $x< x_j$, we have, 
\[
\int_{0}^{1}|E^{j,n}_{l,\theta,\xi}(s)|
                 ds=o(1/\sqrt{n}).
\]
\end{lemme}
\begin{pf}.
In view of Lemma~\ref{lem-bord}, it suffices to show that 
$\int_{\lambda_j/n}^{1-\lambda_j /n}|E^{j,n}_{l,\theta,\xi}(s)|ds=o(1/\sqrt{n})$.
From \eqref{19}, we have, for $n$ large enough,
\begin{eqnarray*}
\int_{\lambda_j/n}^{1-\lambda_j /n}\!\!\!|E^{j,n}_{l,\theta,\xi}(s)|ds&\leq &
  C\,\E\left(e^{|\gamma_j|\left(R_1(\lambda_j/n)+R_2(\lambda_j/n)\right)}\!\!\!\int_0^\infty\!\!\!dm\, 
  \indicatrice_{\{\tilde{x}_j\leq m\leq \tilde{x}_j+R^*(\lambda_j/n)\}}\right)\\
  &\leq &
     C\,\E\left(e^{|\gamma_j|\left(R_1(\lambda_j/n)+R_2(\lambda_j/n)\right)}R^*(\lambda_j/n)
  \indicatrice_{\{ R^*(\lambda_j/n)\geq-\tilde{x}_j\}}\right)\\
  &=&C\sqrt{\frac{\lambda_j}{n}}
      \E\left(e^{|\gamma_j|\sqrt{\frac{\lambda_j}{n}}\left(R_1(1)+R_2(1)\right)}R^*(1)
  \indicatrice_{\{ R^*(1)\geq-\tilde{x}_j\sqrt{n/\lambda_j}\}}\right)
\end{eqnarray*}
Since $\tilde{x}_j<0$, the right hand side of the last equality is $o(1/\sqrt{n})$.
\end{pf}

We will now study the case $x_j<x$. We go back to the decomposition
 \eqref{dec-E} and assume that $n$ is large enough, so that
$\lambda_j J/n<1/4$. Note that, for $s\in[\lambda_j(J+1)/n,1-\lambda_j J/n]$,
we have $N^1_n\geq J$ and $N^2_n\geq J$. So, we have
\[
R^{j,n}_s\leq R^{j,J,n}_s,
\]
where
\[
R^{j,J,n}_s=\min_{-J\leq k\leq J}\check{R}(d^j_n(s)+\lambda_j(k/n)).
\]
%

\begin{lemme}\label{lem-8.3}
If $x_j< x$, we have, 
\[
\lim_{J\to+\infty}\limsup_{n\to\infty}\sqrt{n}\left(
\int_{\frac{\lambda_j (J+1)}{n}}^{1-\frac{\lambda_j J}{n}}\left(E^{j,J,n}_{l,\theta,\xi}(s)-E^{j,n}_{l,\theta,\xi}(s)
                 \right)ds\right)=0,
\]
where
\[
E^{j,J,n}_{l,\theta,\xi}(s)=\E\left(L^j_s
           \alpha_{j}(R_1(s),R_2(1-s))
           \indicatrice_{\left\{\tilde{x}_j\leq  R_1(s)\leq\tilde{x}_j+R^{j,J,n}_s\right\}}\right).
\]
\end{lemme}

\begin{pf}.
Note that
\[
E^{j,J,n}_{l,\theta,\xi}(s)-E^{j,n}_{l,\theta,\xi}(s)=
\E\left(L^j_s
           \alpha_{j}(R_1(s),R_2(1-s))
           \indicatrice_{\Delta^{j,J,n}_s}\right),
\]
with
\begin{eqnarray*}
\Delta^{j,J,n}_s&=&\left\{\tilde{x}_j\leq  R_1(s)\leq\tilde{x}_j+R^{j,J,n}_s\mbox{ and }\exists k\in[0,N^1_n]\cup [-N^2_n,-1],\right.
        \\&&\left.\tilde{x}_j+\check{R}(d^j_n(s)+\lambda_j(k/n))<R_1(s)\right\}.
\end{eqnarray*}
Introducing the notation, for $i=1,2$, and for any real numbers $s_1$, $s_2$ with
$0<s_1<s_2$,
\[
R^\sharp_i(s_1,s_2)=\min_{u\in[s_1,s_2]}R_i(u),
\]
we see that
\begin{eqnarray*}
\Delta^{j,J,n}_s&\subset& \left\{\tilde{x}_j\leq  R_1(s)\leq\tilde{x}_j+R^{j,J,n}_s\right\}
\\&&\cap\left\{R^\sharp_1(\lambda_j J/n,s)<R^{j,J,n}_s \mbox{ or }
    R^\sharp_2(\lambda_j J/n,1-s)<R^{j,J,n}_s\right\}.
\end{eqnarray*}
Note that $R^{j,J,n}_s\leq R^*_1(\lambda_j/n)\wedge R^*_2(\lambda_j/n)$ and 
       $R^\sharp_2(s_1,s_2)\leq R^\sharp_2(s_1)$, where
\[
R^\sharp_2(s)=\min_{u\in[s,+\infty[}R^\sharp_2(u),\quad s\geq 0.
\]
So, we have
\[
|E^{j,J,n}_{l,\theta,\xi}(s)-E^{j,n}_{l,\theta,\xi}(s)|\leq ||F||_\infty 
       \left(F^{j,J,n}_{l,\theta,\xi}(s)+G^{j,J,n}_{l,\theta,\xi}(s)\right),
\]
where
\[
F^{j,J,n}_{l,\theta,\xi}(s)=\E\left(L^j_s
           \indicatrice_{\left\{\tilde{x}_j\leq  R_1(s)\leq\tilde{x}_j+R_1^*(\lambda_j/n),
                R^\sharp_2(\lambda_j J/n)<R_1^*(\lambda_j/n) \right\}}\right)
\]
and
\[
G^{j,J,n}_{l,\theta,\xi}(s)=\E\left(L^j_s
           \indicatrice_{\left\{\tilde{x}_j\leq  R_1(s)\leq\tilde{x}_j+R_2^*(\lambda_j/n),
                R^\sharp_1(\lambda_j J/n,s)<R_2^*(\lambda_j/n) \right\}}\right).
\]
In the sequel we denote by $(\F^i_s)_{s\geq 0}$ ($i=1,2$) the natural filtration of the process 
$(R_i(s))_{s\geq 0}$
and by $\F^i$ the $\sigma$-algebra generated by the union of the 
$\sigma$-algebras $\F^i_s$, $s\geq 0$.

In order to estimate $F^{j,J,n}_{l,\theta,\xi}(s)$, we write
$
L^j_s=\frac{e^{\gamma_jR_1(s)-\frac{\gamma_j^2}{2}(s)}}{R_1(s)}M^j_{1-s},
$
with
\[
M^j_{1-s}=\frac{e^{-\gamma_jR_2(1-s)-\frac{\gamma_j^2}{2}(1-s)}}{2R_2(1-s)}.
\]
By conditioning  with respect to $\F^1$ and using the Cauchy-Schwarz inequality, we get
\begin{eqnarray*}
&&\E\left(\frac{e^{-\gamma_jR_2(1-s)-\frac{\gamma_j^2}{2}(1-s)}}{2R_2(1-s)}
                \indicatrice_{\left\{R^\sharp_2(\lambda_j J/n)<R_1^*(\lambda_j/n) \right\}}\;\vert \;\F^1\right)
                \\&&\leq ||M^j_{1-s}||_2
                \left(\E\left(\indicatrice_{\left\{R^\sharp_2(\lambda_j J/n)<R_1^*(\lambda_j/n) \right\}}
                                                 \;\vert\; \F^1\right)\right)^{1/2}\nonumber\\
                &&=||M^j_{1-s}||_2\sqrt{f_J^{j,n}(R_1^*(\lambda_j/n))},\label{ineg-cond}
\end{eqnarray*}
where $f_J^{j,n}$ is the cumulative distribution function of the random variable 
$R^\sharp_2(\lambda_j J/n)$. Note that, by scaling, we have, for any $r>0$,
 \begin{equation}\label{f_J}
  f_J^{j,n}(r)=f_J\left(r\sqrt{n/\lambda_j}\right),
 \end{equation}
where $f_J$ denotes the cumulative distribution function of $R^\sharp_2(J)$.
On the other hand, we have 
\begin{eqnarray}
||M^j_{1-s}||_2&=&\frac{1}{2\sqrt{1-s}}\left|\left|\frac{e^{-\gamma_j\sqrt{1-s}R_2(1)-\frac{\gamma_j^2}{2}(1-s)}}{R_2(1)}\right|\right|_2
              \nonumber\\
      &=&\frac{1}{2\sqrt{1-s}}
             \left(\int_0^{\infty}e^{-2\gamma_j\sqrt{1-s}m-\gamma_j^2(1-s)-\frac{m^2}{2}}\sqrt{\frac{2}{\pi}}dm
                 \right)^{1/2} \nonumber\\
                 &\leq &
                 \frac{1}{2\sqrt{1-s}}
             \left(\int_{-\infty}^{\infty}e^{-2\gamma_j\sqrt{1-s}m-\gamma_j^2(1-s)-\frac{m^2}{2}}\sqrt{\frac{2}{\pi}}dm
                 \right)^{1/2} \nonumber\\
                 &=&\frac{1}{2\sqrt{1-s}}
             \left(2e^{\gamma_j^2(1-s)}\right)^{1/2}=\frac{1}{\sqrt{2(1-s)}}
            e^{\gamma_j^2(1-s)/2}.\label{Mj}
\end{eqnarray} 
Using the inequalities (\ref{ineg-cond}) and (\ref{Mj}) 
in the expression of $F^{j,J,n}_{l,\theta,\xi}(s)$, we get,
after conditioning with respect to  $\F^1_{\lambda_j/n}$,
\begin{eqnarray*}
F^{j,J,n}_{l,\theta,\xi}(s)&\leq &
  \frac{e^{\gamma_j^2/2}}{\sqrt{2(1-s)}} \E\left(\sqrt{f_J^{j,n}(R_1^*(\lambda_j/n))}\right.
  \\&&\left.\times\int_0^\infty \bar{q}^ j_{s-\lambda_j/n}(R_1(\lambda_j/n),m)
      \indicatrice_{\left\{\tilde{x}_j\leq  m\leq\tilde{x}_j+R_1^*(\lambda_j/n)\right\}}dm\right),
\end{eqnarray*}
with
\[
\bar{q}^ j_{s-\lambda_j/n}(r,m)=e^{\gamma_jm}\bar{q}_{s-\lambda_j/n}(r,m).
\]
We can prove (by the same arguments as in the proof of Proposition~\ref{prop-dom}) that
\[
\int_{\lambda_j(J+1)/n}^{1-\lambda_j J/n}\bar{q}^ j_{s-\lambda_j/n}(r,m)\frac{ds}{\sqrt{2(1-s)}}
  \leq C_{t,\gamma,\sigma}e^{r|\gamma_j|}.
\]
Hence
\begin{eqnarray*}
\int_{\lambda_j(J+1)/n}^{1-\lambda_j J/n}F^{j,J,n}_{l,\theta,\xi}(s) ds
  & \leq &C_{t,\gamma,\sigma}\E\left(\sqrt{f_J^{j,n}(R_1^*(\lambda_j/n))}
  e^{|\gamma_j|R_1(\lambda_j/n)}
      R_1^*(\lambda_j/n)\right)\\
      &=&C_{t,\gamma,\sigma}\E\left(\sqrt{f_J(R_1^*(1))}
  e^{|\gamma_j|\sqrt{\lambda_j/n}R_1(1)}
     \sqrt{\lambda_j/n} R_1^*(1)\right),
\end{eqnarray*}
where the last inequality follows from \eqref{f_J} and the scaling invariance of $R_1$. Hence
\[
\limsup_{n\to \infty}\sqrt{n}\int_{\lambda_j(J+1)/n}^{1-\lambda_j J/n}F^{j,J,n}_{l,\theta,\xi}(s) ds
    \leq C_{t,\gamma,\sigma}\sqrt{\lambda_j}\E\left(\sqrt{f_J(R_1^*(1))}
        R_1^*(1)\right).
\]
We have $f_J(r)=\P(R_2^\sharp(J)\leq r)=\P(R_2^\sharp(1)\leq r/\sqrt{J})$. So, for any $r>0$, \\$\lim_{J\to\infty}f_J(r)=0$. We deduce, by dominated convergence that
\[
\lim_{J\to\infty}\limsup_{n\to \infty}\sqrt{n}\int_{\lambda_j(J+1)/n}^{1-\lambda_j J/n}F^{j,J,n}_{l,\theta,\xi}(s) ds=0.
\]
We will now prove the same property for $G^{j,J,n}_{l,\theta,\xi}(s)$. We have, by conditioning 
with respect to the $\sigma$-algebra generated by
 $\F^2$ and the pair $(R_1(\lambda_j J/n),R_1(s))$,
\begin{eqnarray*}
G^{j,J,n}_{l,\theta,\xi}(s)&=&\E\left(L^j_s
           \indicatrice_{\left\{\tilde{x}_j\leq  R_1(s)\leq\tilde{x}_j+R_2^*(\lambda_j/n),
                R^\sharp_1(\lambda_j J/n,s)<R_2^*(\lambda_j/n) \right\}}\right)\\
                &=&\E\left(L^j_s
           \indicatrice_{\left\{\tilde{x}_j\leq  R_1(s)\leq\tilde{x}_j+R_2^*(\lambda_j/n)\right\}}
                \Delta\left(\frac{\lambda_jJ}{n},s, R_1\left(\frac{\lambda_j J}{n}\right),
                 R_1(s),R_2^*\left(\frac{\lambda_j}{n}\right)\right)\right),
\end{eqnarray*}\normalsize
with the notation, for $0<t_1<t_2$ and $b, y,m>0$,
\[
 \Delta(t_1,t_2, y,m,b)=\P\left(R^\sharp_1(t_1,t_2)\leq b \;|\; R_1(t_1)=y,R_1(t_2)=m\right).
\]
By Lemma~\ref{lem-Bessel}, we have $\Delta(t_1,t_2, y,m,b)\leq \bar\Delta(y,m,b)$,
with
\[
\bar\Delta(y,m,b)=\frac{b(m+y)}{ym}\wedge 1.
\]
Therefore
\begin{eqnarray*}
G^{j,J,n}_{l,\theta,\xi}(s)&\leq&\E\left(L^j_s
           \indicatrice_{\left\{\tilde{x}_j\leq  R_1(s)\leq\tilde{x}_j+R_2^*(\lambda_j/n)\right\}}
                \bar\Delta(R_1(\lambda_j J/n),R_1(s),R_2^*(\lambda_j/n) )\right).
\end{eqnarray*}
Note that
\begin{eqnarray*}
\E\left(\frac{e^{-\gamma_jR_2(1-s)-\frac{\gamma_j^2}{2}(1-s)}}{2R_2(1-s)}\;|\; \F^2_{\lambda_j/n}\vee \F^1\right)
    &=&\frac{1}{2}\int_0^\infty dm e^{-\gamma_j m-\frac{\gamma_j^2}{2}(1-s)}
    \\&&\times\bar{q}_{1-s-\lambda_j/n}(R_2(\lambda_j/n),m)\\
    &\leq &\frac{C}{\sqrt{1-s-\lambda_j/n}}e^{|\gamma_j|R_2(\lambda_j/n)},
\end{eqnarray*}
where we have used  the third estimate of Proposition~\ref{prop-barq}. Now, condition with respect 
to $\F^1_{\lambda_jJ/n}\vee \F^2$ to get
(introducing the random interval $I^{*j}_n=[\tilde{x}_j,\tilde{x}_j+R_2^*(\lambda_j/n)]$ in notations)
\begin{eqnarray*}
G^{j,J,n}_{l,\theta,\xi}(s)&\leq&C\E\left(
      \frac{e^{|\gamma_j|(R_2(\lambda_j/n)+R_1(s))}}{\sqrt{1-s-\lambda_j/n}R_1(s)}
           \indicatrice_{I^{*j}_n}(R_1(s))\right.
                \\&&\left.\times\bar\Delta(R_1(\lambda_j J/n),R_1(s),R_2^*(\lambda_j/n) )\right)\\
                &=&
                C\E\left(
      \frac{e^{|\gamma_j|R_2(\lambda_j/n)}}{\sqrt{1-s-\lambda_j/n}}\int_0^\infty
            \!\!\!\!\!\!dm \,\bar{q}_{s-\lambda_jJ/n}(R_1(\lambda_j J/n),m)\right.
\\&&\left.\times\bar\Delta^j_n(R_1(\lambda_j J/n),m,R_2^*(\lambda_j/n) )\right),              
\end{eqnarray*}
where we have set
\[
\bar\Delta^j_n(r_1,m,r_2 )=e^{|\gamma_j|m}\indicatrice_{I^{*j}_n}(m)\bar\Delta(r_1,m,r_2).
\]
Note that (with the arguments of the proof of Proposition~\ref{prop-dom})
\[
\int_{\lambda_j(J+1)/n}^{1-\lambda_j J/n}e^{|\gamma_j|m}\bar{q}_{s-\lambda_jJ/n}(r,m)\frac{ds}{\sqrt{1-s-\lambda_j/n}}
  \leq C_{t,\gamma,\sigma}e^{r|\gamma_j|},
\]
so that
\begin{eqnarray*}
\int_{\lambda_j(J+1)/n}^{1-\lambda_j J/n} ds G^{j,J,n}_{l,\theta,\xi}(s)&\leq&C_{t,\gamma,\sigma}\E\left(
      e^{|\gamma_j|(R_2(\lambda_j/n)+R_1(\lambda_jJ/n))}
          \!\! \int_0^\infty \!\!\!\! dm\,
            \indicatrice_{I^{*j}_n}(m)\right.
                \\&&\left.\times\bar\Delta(R_1(\lambda_j J/n),m,R_2^*(\lambda_j/n) )\right)\\
                &=&
        C_{t,\gamma,\sigma}\E\left(
      e^{|\gamma_j|\sqrt{\lambda_j/n}(R_1(J)+R_2(1))}
          \!\! \int_0^\infty \!\!\!\! dm\,
            \indicatrice_{\tilde I^{*j}_n}(m)\right.
               \\&&\left.\times \tilde\Delta^j_n(R_1(J),m,R_2^*(1) )\right),                 
\end{eqnarray*}
with
\[
\tilde I^{*j}_n=\left[\tilde{x}_j,\tilde{x}_j+\sqrt{\lambda_j/n}R_2^*(1)\right],
\quad \tilde\Delta^j_n(r_1,m,r_2)=\bar\Delta\left(\sqrt{\lambda_j/n}r_1,m,\sqrt{\lambda_j/n}r_2\right).
\]
We have
\begin{eqnarray*}
                \tilde\Delta^j_n(R_1(J),m,R_2^*(1) )&=&
            \left[R_2^*(1)\left(\frac{1}{R_1(J)}+\frac{\sqrt{\lambda_j/n}}{m}\right)\right]
                \wedge 1    \\
                &\leq &\frac{R_2^*(1)}{R_1(J)}+\left(R_2^*(1) \frac{\sqrt{\lambda_j/n}}{m}\right)\wedge 1.                     
\end{eqnarray*}
By assumption, $x_j< x$, so that  $\tilde{x}_j>0$. Therefore, for $m\in I^{*j}_n$, 
we have $ 1/m\leq 1/\tilde{x}_j$,
   so that
\[
\int_0^\infty \!\!\!\! dm\,
            \indicatrice_{\tilde I^{*j}_n}(m)
                \tilde\Delta^j_n(R_1(J),m,R_2^*(1) )\leq
                \sqrt{\lambda_j/n}R_2^*(1)\left( \frac{R_2^*(1)}{R_1(J)}+R_2^*(1) 
                          \frac{\sqrt{\lambda_j/n}}{\tilde{x}_j}
                \right).
\] 
Hence
\[
\limsup_{n\to\infty}\sqrt{n}\int_{\lambda_j(J+1)/n}^{1-\lambda_j J/n} ds G^{j,J,n}_{l,\theta,\xi}(s)\leq
         C_{t,\gamma,\sigma}\sqrt{\lambda_j}\E\left(\frac{R_2^*(1)^2}{R_1(J)}\right).
\]
Since $\lim_{J\to \infty}R_1(J)=0$, we have
\[
\lim_{J\to \infty}\limsup_{n\to\infty}\sqrt{n}\int_{\lambda_j(J+1)/n}^{1-\lambda_j J/n} ds G^{j,J,n}_{l,\theta,\xi}(s)=0.
\]
\end{pf}

We will now  study the asymptotic behavior of
$\int_{\frac{\lambda_j (J+1)}{n}}^{1-\frac{\lambda_j J}{n}}E^{j,J,n}_{l,\theta,\xi}(s)ds$.
Note that, by conditioning with respect to $\F^1_{\lambda_j  J/n}\vee \F^2_{\lambda_j  J/n}$, we have
\begin{eqnarray*}
E^{j,J,n}_{l,\theta,\xi}(s)&=&
    \E\left(L^j_s
           \alpha_{j}(R_1(s),R_2(1-s))
           \indicatrice_{\left\{\tilde{x}_j\leq  R_1(s)\leq\tilde{x}_j+R^{j,J,n}_s\right\}}\right)\\
           &=&\E\int_0^\infty \!\!\!\!dm\, 
                   \bar{q}^j_{s-\lambda_j J/n}\left(R_1\left({\lambda_jJ/n}\right),m\right)
                      \bar\phi^j_{1-s-\lambda_j J/n}( R_2(\lambda_jJ/n),m)\\&&
                         \times\indicatrice_{\left\{\tilde x_j\leq m\leq \tilde x_j+R^{j,J,n}_s\right\}},
\end{eqnarray*}
with the notations
\[
\bar{q}^j_s(r,m)=e^{\gamma_jm-\frac{\gamma_j^2}{2}}\bar q_s(r,m),\quad
\bar{\phi}^j_s(r,m)=\frac{1}{2}\int_0^\infty\! \bar q_s(r,y)e^{-\gamma_jy}\alpha_j(m,y)dy.
\]

\begin{lemme}
\label{lem-hatE}
Assume $\tilde{x}_j>0$. We have, for any integer $J\geq 1$,
\[
\lim_{n\to\infty}\sqrt{n}\int_{\frac{\lambda_j (J+1)}{n}}^{1-\frac{\lambda_j J}{n}}\left|E^{j,J,n}_{l,\theta,\xi}(s)
                                   -\hat E^{j,J,n}_{l,\theta,\xi}(s)\right|ds=0,
\]
where
\[
\hat E^{j,J,n}_{l,\theta,\xi}(s)=\E\left(\int_0^\infty dm \,\bar q^j_{s}(0,m)
    \bar{\phi}^j_{1-s}(0,m)\indicatrice_{\left\{\tilde x_j\leq m\leq \tilde x_j+R^{j,J,n}_s\right\}}\right).
\]
\end{lemme}

\begin{pf}.
We first  consider
\[
\tilde E^{j,J,n}_{l,\theta,\xi}(s)=\E\left(\int_0^\infty dm \,\bar q^j_{s-\lambda_jJ/n}(0,m)
     \bar{\phi}^j_{1-s-\lambda_jJ/n}(0,m)
      \indicatrice_{\left\{\tilde x_j\leq m\leq \tilde x_j+R^{j,J,n}_s\right\}}\right).
\]
Let $Z^{j,J,n}_{l,\theta,\xi}(s)=E^{j,J,n}_{l,\theta,\xi}(s)
                                   -\tilde E^{j,J,n}_{l,\theta,\xi}(s)$ and
\[
  \zeta^j_s(r_1,r_2,m)=\bar{q}^j_{s-\lambda_jJ/n}(r_1,m)\bar{\phi}^j_{1-s-\lambda_jJ/n}(r_2,m)-
  \bar{q}^j_{s-\lambda_jJ/n}(0,m)\bar{\phi}^j_{1-s-\lambda_jJ/n}(0,m).
\]
We have
\begin{eqnarray*}
Z^{j,J,n}_{l,\theta,\xi}(s)&=&
                                      \E\left(\int_0^\infty \!\!dm 
                                      \indicatrice_{\left\{\tilde x_j\leq m\leq \tilde x_j+R^{j,J,n}_s\right\}}
                                      \zeta^j_s(R_1(\lambda_j J/n),R_2(\lambda_j J/n),m) \right).
\end{eqnarray*}
Now, for all non-negative $r_1$, $r_2$, $m$, 
\begin{eqnarray*}
|\zeta^j_s(r_1,r_2,m)|&\leq&
   |\bar{q}^j_{s-\lambda_jJ/n}(r_1,m)- \bar{q}^j_{s-\lambda_jJ/n}(0,m)|
            |\bar{\phi}^j_{1-s-\lambda_jJ/n}(r_2,m)|\\
   &&+\;
  \bar{q}^j_{s-\lambda_jJ/n}(0,m)|\bar{\phi}^j_{1-s-\lambda_jJ/n}(r_2,m)-\bar{\phi}^j_{1-s-\lambda_jJ/n}(0,m)|.                             
\end{eqnarray*}
By the arguments of Proposition~\ref{prop-barq}, we can easily prove that
\[
\left|\frac{\partial \bar q^j_s}{\partial r}(r,m)\right|\leq \frac{C}{s^{3/2}} e^{r|\gamma_j|},
\quad \left|\frac{\partial \bar \phi^j_s}{\partial r}(r,m)\right|\leq \frac{C}{s} e^{r|\gamma_j|+(\gamma_j^2/2)}.
\]
We deduce that
\begin{eqnarray*}
|\zeta^j_s(r_1,r_2,m)|&\leq&Ce^{(r_1+r_2)|\gamma_j|+(\gamma_j^2/2)}
    r_1\vee r_2\delta^n_j(s),
    \end{eqnarray*}
    with
    \[
    \delta_j^n(s)=\frac{1}{(s-\lambda_jJ/n)^{3/2}(1-s-\lambda_jJ/n)^{1/2}}+
               \frac{1}{(s-\lambda_jJ/n)(1-s-\lambda_jJ/n)}.
    \]
Hence
    \begin{eqnarray*}
\left|Z^{j,J,n}_{l,\theta,\xi}(s)\right|&\leq&C \E\left(R^{j,J,n}_se^{|\gamma_j|(R_1(\lambda_j J/n)+R_2(\lambda_j J/n))+(\gamma_j^2/2)}\right.
   \\&&\left.\times R_1(\lambda_jJ/n)\vee R_2(\lambda_jJ/n) \delta_j^n(s)\right)\\
   &\leq &C_{t,\gamma,\sigma} \E\left(R^*(\lambda_j/n)e^{|\gamma_j|(R_1(\lambda_j J/n)+R_2(\lambda_j J/n))}\right.
   \\&&\left.\times R_1(\lambda_jJ/n)\vee R_2(\lambda_jJ/n) \delta_j^n(s)\right)\\
   &=&C_{t,\gamma,\sigma} \frac{\lambda_j}{n}\E\left(R^*(1)e^{|\gamma_j|\sqrt{\lambda_j/n}(R_1( J)+R_2( J))}
   R_1(J)\vee R_2(J) \delta_j^n(s)\right),
    \end{eqnarray*}
    where the last inequality follows by scaling, and $R^*(u)=R_1^*(u)\vee R_2^*(u)$.
    
Now, fix  $\rho\geq 2$. It can easily be verified that,  for $n$ large enough, we have
\begin{eqnarray*}
\int_{\frac{\lambda_j \rho J}{n}}^{1-\frac{\lambda_j\rho J}{n}} \delta_j^n(s)ds
     &\leq & C\left(\frac{\sqrt{n}}{\sqrt{\lambda_j(\rho-1)J}}+\log\left(\frac{n}{\lambda_j(\rho-1) J}\right)\right),
\end{eqnarray*}   
so that
\[
\limsup_{n\to\infty}\sqrt{n}\int_{\frac{\lambda_j \rho J}{n}}^{1-\frac{\lambda_j\rho J}{n}}\left|Z^{j,J,n}_{l,\theta,\xi}(s)\right|ds
\leq \frac{C_{t,\gamma,\sigma}}{\sqrt{\lambda_j(\rho-1)J}}.
\] 
On the other hand, it can be proved (as in Lemma~\ref{lem-bord}) that, for any fixed $\rho\geq 2$,
\[
\int_{\frac{\lambda_j (J+1)}{n}}^{\frac{\lambda_j \rho J}{n}}\left|Z^{j,J,n}_{l,\theta,\xi}(s)\right|ds
+\int_{1-\frac{\lambda_j\rho J}{n}}^{1-\frac{\lambda_jJ}{n}}\left|Z^{j,J,n}_{l,\theta,\xi}(s)\right|ds=o(1/\sqrt{n}).
\]
Therefore
\[
\limsup_{n\to\infty}\sqrt{n}\int_{\frac{\lambda_j (J+1)}{n}}^{1-\frac{\lambda_jJ}{n}}\left|Z^{j,J,n}_{l,\theta,\xi}(s)\right|ds
\leq \frac{C_{t,\gamma,\sigma}}{\sqrt{\lambda_j(\rho-1)J}},
\]
and, by letting $\rho$ go to infinity, we conclude that
\[
\lim_{n\to\infty}\sqrt{n}\int_{\frac{\lambda_j (J+1)}{n}}^{1-\frac{\lambda_j J}{n}}\left|E^{j,J,n}_{l,\theta,\xi}(s)
                                   -\tilde E^{j,J,n}_{l,\theta,\xi}(s)\right|ds=0.
\] 
It remains to show that
\[
\lim_{n\to\infty}\sqrt{n}\int_{\frac{\lambda_j (J+1)}{n}}^{1-\frac{\lambda_j J}{n}}\left|\hat E^{j,J,n}_{l,\theta,\xi}(s)
                                   -\tilde E^{j,J,n}_{l,\theta,\xi}(s)\right|ds=0.
\]
 We have, for 
$s\in\left(\frac{\lambda_j (J+1)}{n},1-\frac{\lambda_j J}{n}\right)$,
\begin{eqnarray*}
\left|\hat E^{j,J,n}_{l,\theta,\xi}(s)
                                   -\tilde E^{j,J,n}_{l,\theta,\xi}(s)\right|
                                 &=&
                                   \left|\E\left(\int_{\tilde x_j}^{ \tilde x_j+R^{j,J,n}_s}\!\!\!\! dm\,\eta^n_j(s,m)
    \right)\right|\leq
       \E\left(\int_{\tilde x_j}^{ \tilde x_j+\frac{R^*(\lambda_j)}{\sqrt{n}}}\!\!\!\! dm\,|\eta^n_j(s,m)|
    \right),
\end{eqnarray*}
with
\begin{eqnarray*}
\eta^n_j(s,m)&=&
\bar q^j_{s}(0,m)\bar{\phi}^j_{1-s}(0,m)
                                   -\bar q^j_{s-\lambda_jJ/n}(0,m)\bar{\phi}^j_{1-s-\lambda_jJ/n}(0,m).
\end{eqnarray*}
Recall that $\bar q_s(0,m)=\sqrt{2/\pi}\frac{m}{s^{3/2}}e^{-m^2/2s}$. 
Note that, if $\tilde x_j\leq m\leq \tilde x_j+\frac{R^*(\lambda_j)}{\sqrt{n}}$, we have
 $e^{-m^2/4s}\leq e^{-\tilde{x}_j^2/4s}$, so that, for $s\in(0,1)$,
\begin{eqnarray}\label{29}
\bar q^j_{s}(0,m)\leq \frac{e^{-\tilde{x}_j^2/4s}}{s^{3/2}}me^{|\gamma_j|m-(m^2/4)}\leq C_j,
\end{eqnarray}
for some positive constant $C_j$ (depending on $\tilde x_j$, but not on $s$ or $m$).
 
  Furthermore
\begin{eqnarray*}
\bar{\phi}^j_s(0,m)&=&\frac{1}{2}\int_0^\infty dy\bar q_s(0,y)e^{-\gamma_j y}\alpha_j(m,y)\\
             &=&\frac{1}{2\sqrt{s}}\int_0^\infty dy\sqrt{2/\pi}ye^{-y^2/2-\gamma_j\sqrt{s}y}
                    \alpha_j(m,\sqrt{s}y),
             \end{eqnarray*}
so that $|\bar{\phi}^j_s(0,m)|\leq C/\sqrt{s}$, for some $C>0$.

We deduce thereof that, for 
 $s\in\left(\frac{\lambda_j (J+1)}{n},1-\frac{\lambda_j J}{n}\right)$,
 \[
 \int_{\tilde x_j}^{ \tilde x_j+\frac{R^*(\lambda_j)}{\sqrt{n}}}\!\!\!\!\!\!\!dm\,|\eta^n_j(s,m)|
    \leq C_j\frac{R^*(\lambda_j)}{\sqrt{n}}
                                  \frac{1}{\sqrt{1-s-\lambda_jJ/n}}.
 \]
 From this estimate, together with the fact that for a fixed $s$ we have 
 $\int_{\tilde x_j}^{ \tilde x_j+\frac{R^*(\lambda_j)}{\sqrt{n}}}\!\!\!\!\!\!\!dm\,|\eta^n_j(s,m)|=o(1/\sqrt{n})$
 a.s., we easily deduce that 
 \[
\lim_{n\to\infty}\sqrt{n}\int_{\frac{\lambda_j (J+1)}{n}}^{1-\frac{\lambda_j J}{n}}
\left|\eta^n_j(s,m)\right|ds=0.
\]
 \end{pf}

\begin{lemme}\label{lem-fin}
We have, if $\tilde{x}_j>0$,
\[
\int_{\frac{\lambda_j (J+1)}{n}}^{1-\frac{\lambda_j J}{n}}\hat E^{j,J,n}_{l,\theta,\xi}(s)ds=
\sqrt{\lambda_j/n}\E(R^J)\int_0^1\varphi^j_s(\tilde{x}_j)ds+o(1/\sqrt{n}),
\]
where
\[
\varphi^j_s(m)=\bar q^j_s(0,m)\bar{\phi}^j_{1-s}(0,m),\quad R^J=\min_{|k|\leq J}\check{R}(U+k),
\]
and the random variable  $U$ is uniformly distributed  on $[0,1]$ and independent of $\check{R}$.
\end{lemme}

\begin{pf}.
We have
\begin{eqnarray*}
\hat E^{j,J,n}_{l,\theta,\xi}(s)=\E \int_{\tilde x_j}^{\tilde x_j+R^{j,J,n}_s}  \!\!\!\!dm\,\varphi^j_s(m),
\end{eqnarray*}
so that, intoducing the notation
\[
\bar{E}^{j,J,n}_{l,\theta,\xi}(s)=\E(R^{j,J,n}_s)\bar q^j_s(0,\tilde{x}_j)\bar{\phi}^j_{1-s}(0,\tilde x_j)=
      \E(R^{j,J,n}_s)\varphi^j_s(\tilde{x}_j),
\]
we can write, for $s\in(\lambda_j(J+1)/n,1-\lambda_jJ/n)$,
\begin{eqnarray*}
\left|\hat E^{j,J,n}_{l,\theta,\xi}(s)-\bar{E}^{j,J,n}_{l,\theta,\xi}(s)\right|
&\leq&\E \int_{\tilde x_j}^{\tilde x_j+R^{j,J,n}_s}  \!\!\!\!dm\,\left|\varphi^j_s(m)-\varphi^j_s(\tilde{x}_j)\right|
     \\
     &\leq &
     \E \int_{\tilde x_j}^{\tilde x_j+\frac{R^*(\lambda_j)}{\sqrt{n}}}  \!\!\!\!dm
                      \,\left|\varphi^j_s(m)-\varphi^j_s(\tilde{x}_j)\right|.
                      \end{eqnarray*}
Note that, using \eqref{29} and the estimate $|\bar{\phi}^j_s(0,m)|\leq C/\sqrt{s}$, we have,
 for $\tilde x_j\leq m\leq \tilde x_j+R^*(\lambda_j)/\sqrt{n}$,
\[
\left|\varphi^j_s(m)\right|\leq \frac{C_j}{\sqrt{1-s}}.
\]
Now, for a fixed $s\in (0,1)$, $m\mapsto \varphi^j_s(m)$ is continuous.
Indeed the continuity of $\bar{q}^j_s(0,.)$ is obvious and the continuity of $\bar{\phi}^j_s(0,.)$
follows from the continuity assumption on $F$ and the equalities
\begin{eqnarray*}
\bar{\phi}^j_s(0,m)&=&\frac{1}{2}\int_0^\infty dy\bar q_s(0,y)e^{-\gamma_j y}\alpha_j(m,y)\\
            &=&\frac{1}{2}\int_0^\infty dy\bar q_s(0,y)e^{-\gamma_j y}F(x_j+\sigma_j(m-y),x_j+\sigma_jm)\\
            &=&
            \frac{1}{2}\int_{-\infty}^m dz\bar q_s(0,m-z)e^{-\gamma_j (m-z)}F(x_j+\sigma_jz,x_j+\sigma_jm).
\end{eqnarray*}
Due to the continuity of $m\mapsto \bar{\phi}^j_s(m)$ we have, for a fixed $s\in (0,1)$,
\[
\int_{\tilde x_j}^{\tilde x_j+\frac{R^*(\lambda_j)}{\sqrt{n}}}  \!\!\!\!dm
                      \,\left|\varphi^j_s(m)-\varphi^j_s(\tilde{x}_j)\right|=o(1/\sqrt{n}) \mbox{ a.s.,}
\]
Hence
\[
\int_{\frac{\lambda_j (J+1)}{n}}^{1-\frac{\lambda_j J}{n}}\hat E^{j,J,n}_{l,\theta,\xi}(s)ds=
\int_{\frac{\lambda_j (J+1)}{n}}^{1-\frac{\lambda_j J}{n}}
\E(R^{j,J,n}_s)\varphi^j_s(\tilde{x}_j)ds+o(1/\sqrt{n}).
\]
We have
\begin{eqnarray*}
\E(R^{j,J,n}_s)&=&\E\left(\min_{|k|\leq J} \check{R}(d_n^j(s)+\lambda_j k/n)\right)\\
             &=&\sqrt{\lambda_j/n}\E\left(\min_{|k|\leq J} \check{R}(nd_n^j(s)/\lambda_j+k)\right)\\
             &=&\sqrt{\lambda_j/n}f\left(\frac{nd_n^j(s)}{\lambda_j}\right),
\end{eqnarray*}
where, for $u\in[0,1]$,  $f(u)=\E\left(\min_{|k|\leq J} \check{R}(u+k)\right)$. 
Using the definition of $d^j_n(s)$ and classical arguments, one can show that, 
for any integrable function $g$ on $[0,1]$, we have
\[
\lim_{n\to\infty}\int_{\frac{\lambda_j (J+1)}{n}}^{1-\frac{\lambda_j J}{n}}
     f\left(\frac{nd_n^j(s)}{\lambda_j}\right)g(s)ds=\int_0^1g(s)ds\int_0^1 f(u)du.
     \]
Therefore
\[
\int_{\frac{\lambda_j (J+1)}{n}}^{1-\frac{\lambda_j J}{n}}\hat E^{j,J,n}_{l,\theta,\xi}(s)ds
    =\sqrt{\lambda_j/n}\E(R^J)\int_0^1 ds \varphi^j_s(\tilde{x}_j)+ o(1/\sqrt{n}),
\]
which proves the Lemma.
\end{pf}

\begin{pf} \itshape of Theorem~\ref{thm-convcond}. \upshape
Recall that
\begin{eqnarray}
\E_{l,\theta,\xi}\left(F(X_{t_{j+1}^-},M^j)
  \indicatrice_{\{ M^j\geq x>M^{j,n}\}}\right)
   &=&
\int_0^1E^{j,n}_{l,\theta,\xi}(s)ds\nonumber\\
&=&\int_{\frac{\lambda_j (J+1)}{n}}^{1-\frac{\lambda_j J}{n}}  E^{j,n}_{l,\theta,\xi}(s)ds
      +o(1/\sqrt{n}),\label{32*}
\end{eqnarray}
for any positive integer $J$. Here, we have used  Lemma~\ref{lem-bord}.
Note that, if $x<x_j$, we have, due to Lemma~\ref{lem-8.3-},
$\E_{l,\theta,\xi}\left(F(X_{t_{j+1}^-},M^j)
  \indicatrice_{\{ M^j\geq x>M^{j,n}\}}\right)=o(1/\sqrt{n})$.

We now assume $x>x_j$ and fix a positive integer $J$. We have, using Lemma~\ref{lem-hatE} and Lemma~\ref{lem-fin},
\begin{eqnarray*}
\int_{\frac{\lambda_j (J+1)}{n}}^{1-\frac{\lambda_j J}{n}}  E^{j,J,n}_{l,\theta,\xi}(s)ds&=&
\int_{\frac{\lambda_j (J+1)}{n}}^{1-\frac{\lambda_j J}{n}}\hat E^{j,J,n}_{l,\theta,\xi}(s)ds+o(1/\sqrt{n})\\
&=&
\sqrt{\frac{\lambda_j}{n}}\E(R^J)\int_0^1 \varphi^j_s(\tilde{x}_j)ds+o(1/\sqrt{n}).
\end{eqnarray*}
Note that $\lim_{J\to \infty}\E(R^J)=\beta_1$, so that
\begin{eqnarray}\label{32}
\lim_{J\to+\infty}\limsup_{n\to\infty}\sqrt{n}\left(
\int_{\frac{\lambda_j (J+1)}{n}}^{1-\frac{\lambda_j J}{n}}{E}^{j,J,n}_{l,\theta,\xi}(s)ds
  -\sqrt{\frac{\lambda_j}{n}}\beta_1\int_0^1 \varphi^j_s(\tilde{x}_j)ds\right)=0.
\end{eqnarray}
By combining \eqref{32*}, \eqref{32} 
and Lemma~\ref{lem-8.3}, we derive
\begin{eqnarray}\label{dl}
\E_{l,\theta,\xi}\left(F(X_{t_{j+1}^-},M^j)
  \indicatrice_{\{ M^j\geq x>M^{j,n}\}}\right)
   &=&
\sqrt{\frac{\lambda_j}{n}}\beta_1\int_0^1 \varphi^j_s(\tilde{x}_j)ds
      +o(1/\sqrt{n}).
\end{eqnarray}
On the other hand, for any $\rho>0$, we have (using Proposition~\ref{prop-rep-Bessel}
for a function which does not depend on the difference $M^j-M^{j,n}$)
\begin{eqnarray*}
\E_{l,\theta,\xi}\left(F(X_{t_{j+1}^-},M^j)
  \indicatrice_{\{ M^j\geq x>M^j-\rho\}}\right)
   &=&
\int_0^1ds\int_0^{\infty}dm\indicatrice_{\{\tilde{x}_j\leq m\leq \tilde{x}_j+\frac{\rho}{\sigma_j}\}} dm\varphi^j_s(m).
 \end{eqnarray*}
 If $\tilde{x}_j<0$, we get
 $\E_{l,\theta,\xi}\left(|F(X_{t_{j+1}^-},M^j)|
  \indicatrice_{\{ M^j\geq x>M^j-\rho\}}\right)=0$ for $\rho<\sigma_j|\tilde{x}_j|$,
  so that 
  \begin{eqnarray*}
\E_{l,\theta,\xi}\left(F(X_{t_{j+1}^-},M^j)
  \indicatrice_{\{ M^j\geq x>M^j-\sigma\beta_1\sqrt{t/n}\}}\right)=o(1/\sqrt{n}).
\end{eqnarray*}

If $\tilde{x}_j>0$, we have
 \begin{eqnarray*}
\E_{l,\theta,\xi}\left(|F(X_{t_{j+1}^-},M^j)|
  \indicatrice_{\{ M^j\geq x>M^j-\rho\}}\right)
&=&
\frac{\rho}{\sigma_j}\int_0^1\varphi^j_s(\tilde{x}_j)ds+o(\rho)
 \end{eqnarray*}
as $\rho$ goes to $0$. Therefore
\begin{eqnarray*}
\E_{l,\theta,\xi}\left(F(X_{t_{j+1}^-},M^j)
  \indicatrice_{\{ M^j\geq x>M^j-\sigma\beta_1\sqrt{t/n}\}}\right)
   &=&
\frac{\sigma\beta_1\sqrt{t}}{\sigma_j\sqrt{n}}\int_0^1\varphi^j_s(\tilde{x}_j)ds
      +o(1/\sqrt{n}).
\\      &=&\frac{\beta_1\sqrt{t}}{\sqrt{t_{j+1}-t_j}\sqrt{n}}
      \int_0^1\varphi^j_s(\tilde{x}_j)ds
      +o(1/\sqrt{n})\\
      &=&
      \beta_1\sqrt{\frac{\lambda_j}{n}}\int_0^1 \varphi^j_s(\tilde{x}_j)ds
      +o(1/\sqrt{n}),
\end{eqnarray*}
which completes the proof of the second statement of the Theorem. The first one
can be proved by the same method.
\end{pf}
\begin{rmq}\label{rem:fin}\rm
It can be deduced from~\eqref{dl} that we have an expansion
\[
\E\left(g(X_t)\indicatrice_{\{M_t\geq x>M^n_t\}}\right)=\frac{C}{\sqrt{n}}+o(1/\sqrt{n}),
\]
for some constant $C$. This can be used to derive an expansion for the difference
between continuous and discrete barrier option prices.
\end{rmq}


\begin{thebibliography}{2}

\bibitem{asmussen-al95} \textsc{Asmussen, S., Glynn P., Pitman J.}: Discretization error in simulation of one-dimensional reflecting brownian motion. {\em The Annals of Applied Probability} {\bf 3}(4), 875-896, (1995).

\bibitem{broadie-glasserman} \textsc{Broadie, M., Glasserman P., Kou S. G.}: A Continuity Correction For Discrete Barrier Options. {\em Mathematical Finance} {\bf 7}(4), 325-348, October (1997).

\bibitem{broadie-glasserman1999} \textsc{Broadie, M., Glasserman P., Kou S. G.}: Connecting
 Discrete and Continuous Path-Dependent Options.  {\em Finance and Stochastics}
  {\bf  3}, 55-82,  (1999).

\bibitem{dia} \textsc{Dia, E. H. A.}:  Exotic Options under Exponential L\'evy Model. Doctoral thesis, 
Universit\'e Paris-Est, http://tel.archives-ouvertes.fr/INSMI/tel-00520583/fr/, (2010).

\bibitem{DiaLamberton} \textsc{Dia, E. H. A., D. Lamberton}: 
Connecting Discrete and Continuous Lookback or Hindsight 
Options in Exponential L\'evy Models, {\em submitted for publication}, 2010. 
 

\bibitem{karatzas-shreve} \textsc{Karatzas, I., Shreve S.E.}: Brownian motion and stochastic calculus, 2nd edition Springer (1991).

\bibitem{kou03} \textsc{Kou, S. G.}: On pricing of discrete barrier options. 
{\em Statistica Sinica} {\bf 13}, 955-964, (2003).

\bibitem{kou02} \textsc{Kou, S. G.}: A jump-diffusion model for option pricing. 
{\em Management Science} {\bf 48} (8), 1086-1101, (2002).

\bibitem{kou-wang04} \textsc{Kou, S. G., Wang H.}: Option pricing under a double exponential jump diffusion model. {\em Management Science} {\bf 50} (9), 1178-1192, (2004).

\bibitem{kyprianou} \textsc{Kyprianou, A. E.}: 
{\em Indroductory Lectures on Fluctuations of L\'evy Processes with Applications}.  Springer-Verlag,
 (2006).

\bibitem{revuz-yor} \textsc{Revuz, D., Yor M.}: 
{\em Continuous martingales and Brownian motion, 2nd edition},
Springer, (1994).

\bibitem{sato} \textsc{Sato, K.}: 
{\em L\'evy processes and infinitely divisible distributions.} Cambridge University Press, (2005).
\end{thebibliography}
\end{document}